\documentclass[11pt]{article}
\usepackage{amsfonts,amssymb,amsmath,numbysec,bbm}
\hsize \textwidth
\advance \hsize by -\marginparwidth
\evensidemargin \oddsidemargin
\usepackage{amssymb}
\usepackage{amsthm}
\usepackage[shortlabels]{enumitem}
	\setlist{itemindent=-1pt, itemsep=-2pt,topsep=2pt}
\usepackage[normalem]{ulem}
\usepackage{cancel}
\usepackage{scalerel}
\usepackage{multicol}
\usepackage{xfrac}
\usepackage{xcolor}
\usepackage{tikz}
	\usetikzlibrary{math}

\usepackage{mathrsfs}

\advance\hoffset by 5mm

\newcommand{\Rm}{{\mathbb R}}
\newcommand{\R}{\mathbb R}
\newcommand{\eps}{\varepsilon}
\newcommand{\commentout}[1]{}

\renewcommand{\phi}{\varphi}

\DeclareMathOperator{\supp}{supp}

\newcommand{\Rd}{R}

\usepackage[colorlinks=true,pdfpagemode=UseNone,urlcolor=blue,linkcolor=blue,citecolor=blue,breaklinks=true]{hyperref}
\usepackage{cleveref}

\newtheorem{thm}{Theorem}[section]
\newtheorem{lem}[thm]{Lemma}

\newtheorem{prop}[thm]{Proposition}

\Crefname{lem}{Lemma}{Lemmas}
\crefname{lem}{Lemma}{Lemmas}
\Crefname{thm}{Theorem}{Theorems}

\newcommand{\farc}{\frac}
\newcommand{\be}{\begin{equation}}
\newcommand{\ee}{\end{equation}}

\newcommand{\bal}{\begin{aligned}}
\newcommand{\enbal}{\end{aligned}}

\newcommand{\one}{{\mathbbm{1}}}
\newcommand{\1}{\one}
\newcommand{\ifnty}{\infty}

\newcommand{\bhat}[1]{\expandafter\hat#1} 

\numberbysection
\numberwithin{equation}{section}

\usepackage{autonum}

\begin{document}

\author{Jing An \and Christopher Henderson \and Lenya Ryzhik}
\title{Front location determines convergence rate to traveling waves}

\maketitle

\begin{abstract} 
We propose a novel method for establishing the convergence rates of solutions to reaction-diffusion equations to traveling waves.  The analysis is based on the study of 
the traveling wave shape defect function introduced in~\cite{AHR2}.  It turns out that the convergence rate is controlled by the distance between the ``phantom front location'' 
for the shape defect function and the true front location of the solution.  Curiously, the convergence to a traveling wave has a pulled nature, regardless of whether the traveling 
wave itself is of pushed, pulled, or pushmi-pullyu type.  In addition to providing new results, this approach simplifies dramatically the proof in the Fisher-KPP case and gives a unified, 
succinct explanation for the known algebraic rates of convergence in the Fisher-KPP case and the exponential rates in the pushed case. 
\end{abstract}

\section{Introduction}\label{s.intro}

We consider the long-time behavior of solutions to reaction-diffusion equations of the form
\be\label{e.rde}
u_t = u_{xx} + f(u),~~t>0,~~x\in\Rm,
\ee
with a nonlinearity $f\in C^2([0,1])$ that satisfies
\be\label{e.f}
	f(0) = f(1) = 0,
	\quad f'(0) >0,
	\quad f(u) > 0 \text{ for } u \in (0,1).
\ee
In addition, we normalize the nonlinearity so that 
\be\label{e.normalization}
f'(0)=1.
\ee
This condition can be achieved by a simple space-time rescaling and is not  
an extra assumption on $f(u)$. 
Reaction-diffusion equations of the form~\eqref{e.rde} are used in a wide variety of settings 
to understand how the interplay of diffusive spreading  and growth gives rise to front propagation and invasions.
Our interest is in precisely quantifying this behavior.

\subsection*{Convergence in shape to a traveling wave}
 
Traveling waves are solutions to (\ref{e.rde}) of the form $u(t,x)=U_c(x-ct)$, with a profile $U_c(x)$ such that 
\be\label{mar2612}
-cU_c'=U_c''+f(U_c),
\ee
and
\be\label{jun2002}
0<U_c(x)<1~~\hbox{for all $x\in\Rm$}, ~~U_c(-\infty)=1,~~U_c(+\ifnty)=0.
\ee
Solutions to~\eqref{mar2612}-\eqref{jun2002} are only unique up to translation, so we often
fix the choice of the wave by
the normalization
\be\label{jun2004}
U_c(0)=\frac12.
\ee
Another natural normalization is mentioned in Section~\ref{s.results}, 
see (\ref{july1002}) below.
For nonlinearities satisfying~(\ref{e.f}), there exists a minimal speed $c_*>0$ such that traveling waves exist if and only if $c\ge c_*$~\cite{HadelerRothe}.
The normalization (\ref{e.normalization}) implies
that $c_*\ge 2$. We denote the profile of the wave corresponding to the minimal front speed $c_*$ as $U_*(x)$. 

The study of the long time behavior of the solutions to (\ref{e.rde}) with initial conditions that decay rapidly as $x\to+\infty$  goes back to
the original papers~\cite{Fisher,kpp}. To be concrete and avoid some additional technicalities, we  momentarily consider the case where the initial
condition for (\ref{e.rde}) is a step-function:
\be\label{mar2302}
u_0(x)=u(0,x)=\one(x\le 0).
\ee
It is well known that this assumption may be greatly relaxed, as long as $u_0(x)$ is sufficiently rapidly decaying as $x\to+\infty$,
see~\cite{berestycki2017exact,BD-2015} for a recent detailed analysis of this issue. 
It was shown in the original KPP paper~\cite{kpp}  that the solution $u(t,x)$ to (\ref{e.rde}) converges to $U_*(x)$
in shape. That is, there exists a reference frame $m(t)$ such that 
\be\label{e.c021601}
u(t,x + m(t)) - U_*(x)=o(1),~~\hbox{as $t\to+\infty$.}
\ee
We will refer to $m(t)$ as the front location. Note that, strictly speaking,
it is only defined up to an~$o(1)$ term as $t\to+\infty$. 
Moreover, the KPP paper showed that the front location $m(t)$ has the asymptotics
\be\label{mar2306}
m(t) = c_*t + o(t),~~\hbox{as $t\to+\infty$.}
\ee
The extraordinarily 
innovative proof in~\cite{kpp} 
relies on, in modern terminology, an intersection number argument 
and 
can be extended not only to all Lipschitz $f(u)$ that satisfy (\ref{e.f}),
but to a much larger classes of nonlinearities. In that sense, both  (\ref{e.c021601}) and (\ref{mar2306}) are fairly universal results. 

\subsubsection*{Front location and convergence rates in the pushed and pulled cases}

On the other hand, both
the precise character of the~$o(t)$ correction to the front location 
in (\ref{mar2306}) and  the rate of  the ``convergence in shape'' in (\ref{e.c021601})
depend heavily on the profile of the nonlinearity~$f(u)$, as neither 
can  be easily obtained from the intersection number arguments.   

The results quantifying these convergence rates and making the asymptotics of the front location~$m(t)$ more precise than (\ref{mar2306})
are more modern and
are very different in what are known as
the ``pushed" and ``pulled" regimes. Recall that,
informally, front propagation is pushed if it is ``bulk dominated"
and is pulled if it is ``tail dominated". For positive nonlinearities that satisfy (\ref{e.f})-(\ref{e.normalization})
the spreading speed for the linearized problem
\be\label{jun2210}
u_t=u_{xx}+u,
\ee
is $c_{\rm lin}=2$.  We will give a more refined definition below but for the moment the reader can think that propagation is pushed
if $c_*>c_{\rm lin}=2$ and pulled if $c_*=c_{\rm lin}=2$.
Contemporary arguments to establish convergence rates in the pushed case are spectral in nature, while, for pulled fronts, are motivated in great part by the connection to branching Brownian motion and other log-correlated random fields,
and typically use entirely different techniques.   

When the front is pushed, so that $c_*>2$, its 
location has the asymptotics
\be\label{mar2308}
m(t)=c_*t+x_0+o(1),~~\hbox{as $t\to+\infty$},
\ee
with some $x_0\in\Rm$. Moreover, the convergence rate in (\ref{e.c021601})
is exponential~\cite{fife1977approach,Rothe}:
\be\label{mar2312}
|u(t,x + m(t)) - U_*(x)| \le ce^{-\omega t},
\ee
with some $\omega>0$. 
The proofs of (\ref{mar2308})-(\ref{mar2312}) in~\cite{fife1977approach,Rothe} as well as the later extensions to other ``pushed fronts" problems
are  based on spectral gap arguments and provide implicit estimates on the exponential rate $\omega>0$ of convergence in~(\ref{mar2312}).

On the other hand, when  $f(u)$ is of the Fisher-KPP type, so that, in addition to (\ref{e.f}), it satisfies 
 \be\label{mar2310}
 f(u)\le f'(0)u,~~\hbox{for all $0<u<1$}, 
 \ee
the propagation is pulled and spreading is dominated by the region
far ahead of the front.  Under this assumption, when the normalization (\ref{e.normalization}) is adopted, the minimal speed~$c_*=c_{\rm lin} = 2$ and the front location has the asymptotics
 \be\label{mar2404}
 m(t)=2t-\farc{3}{2}\log t+x_0+o(1),~~\hbox{as $t\to+\infty$},
 \ee
with some $x_0\in\Rm$, first established in the pioneering works by Bramson~\cite{Bramson1,Bramson2} via the connection with branching Brownian motion.   The Bramson asymptotics
was re-visited in~\cite{AHR,AHR2,Avery,AS2,Giletti,HNRR,Lau,NRR1,Roberts,Uchiyama}, including in some more general pulled settings, and also refined in~\cite{berestycki2017exact,berestycki2018new,Giletti,Graham,NRR2}.  
However, unlike in the pushed case, where the front location asymptotics (\ref{mar2308}) 
was sufficient for the convergence rate estimate (\ref{mar2312}), 
obtaining a convergence rate in (\ref{e.c021601}) for the Fisher-KPP nonlinearities  required a much finer asymptotics than given by the Bramson result~(\ref{mar2404}).  
To this end,
Graham has improved in~\cite{Graham}  the Bramson asymptotics for the Fisher-KPP nonlinearities 
to show that 
\be\label{mar2408}
m(t)=2t-\farc{3}{2}\log t+x_0-\farc{3\sqrt{\pi}}{\sqrt{t}}+\farc{9}{8}(5-6\log 2)\farc{\log t}{t}+\farc{x_1}{t}+o\Big(\farc{1}{t}\Big),~~\hbox{as $t\to+\infty$},
\ee
with some $x_0,x_1\in\Rm$. This confirmed a series of formal predictions in~\cite{berestycki2018new,Ebert-vanSaarlos}, partly proved in~\cite{henderson2016,NRR2}.
The ``very fine'' asymptotics in (\ref{mar2408})  leads to a convergence
bound of the form 
\be\label{mar2410}
 |u(t,x + m(t)) - U_*(x)| =O\Big(\farc{1}{t}\Big) 
\ee
after using an asymptotic expansion based on (\ref{mar2408})
that approximately solves~\eqref{e.rde}.  It was also shown 
in~in~\cite{Graham} that this rate can not be improved
for the Fisher-KPP nonlinearities.  We note that, with different assumptions on the initial data that rule out~\eqref{mar2302} and its compact perturbations, faster convergence rates were proven by Avery and Scheel~\cite{AS1}.

While the Bramson asymptotics (\ref{mar2404}) holds for all Fisher-KPP reactions, it does not hold for all nonlinearities that satisfy (\ref{e.f})-(\ref{e.normalization})
for which $c_*=2$. As was shown in~\cite{AHR2,Giletti}, there is a class of nonlinearities $f(u)$ such that the front location asymptotics is not (\ref{mar2404}) but 
\be\label{jun2208}
m(t) = 2t - \frac{1}{2}\log t + x_0 +o(1),~~\hbox{as $t\to+\infty$}. 
\ee
Informally, this happens when $f(u)$ is exactly at the pushed-pulled transition.  We refer to these as ``pushmi-pullyu" fronts. 
Thus, the distinction between various regimes of propagation can not be made based solely on whenever the propagation speed is predicted by the linearization
(\ref{jun2210}) or not. It turns out that it should be made based both on the propagation speed and the asymptotics behavior of the 
traveling wave as $x\to+\infty$.
Let us, therefore,  define terminology for the three classes roughly discussed above. We remind the reader that $f(u)$ satisfies (\ref{e.f})-(\ref{e.normalization}).
\begin{itemize}[itemindent=-7pt]
	\item A traveling wave is {\bf pushed} if $c_* > 2$.  
	\item A traveling wave is {\bf pulled} if $c_* = 2$ and there is some $A_0 > 0$ such that
		\be\label{e.c062202}
			U_*(x)
				= A_0 x e^{-x} + O(e^{-x})
				\qquad\text{ as } x \to \infty.
		\ee
	\item A traveling wave is {\bf pushmi-pullyu} if $c_* = 2$ and there is $A_1>0$ such that
		\be\label{e.c062203}
			U_*(x)
				= A_1 e^{-x} + o(e^{-x})
				\qquad\text{ as } x \to \infty.
		\ee
\end{itemize}
We refer the reader to~\cite{AHR2,AHS,berestycki2017exact,berestycki2018new,Ebert-vanSaarlos,garnier2012inside,Giletti,vanSaarlos} for more in depth discussion.    We often abuse terminology and refer to the nonlinearity itself
as being ``pushed,'' ``pulled,'' or ``pushmi-pullyu.''

A simple linearization argument shows that the two asymptotics in~\eqref{e.c062202}-\eqref{e.c062203} are the only possibilities when $c_* = 2$, so the cases above are exhaustive.  Intuitively,  once the normalization~(\ref{e.normalization}) is fixed,
``large'' nonlinearities $f$ correspond to pushed fronts, ``small'' ones correspond to pulled
fronts, and the boundary case corresponds to pushmi-pullyu fronts.

There are two important points to make before discussing our results.  First, while convergence rates have been established in the Fisher-KPP and pushed cases, nothing quantitative is known for the intermediate cases; that is, pushmi-pullyu nonlinearities  and pulled nonlinearities not satisfying the Fisher-KPP condition~\eqref{mar2310}.  Second, the arguments used to establish convergence rates in the Fisher-KPP and pushed regimes are quite different.  This indicates the difficulty in closing the gap: establishing sharp rates in the transitional cases and developing a cohesive understanding of convergence rates in all cases.

\subsubsection*{An informal statement of the results}

Our interest here is to complete and unify the separate pictures for the pulled, pushed, and pushmi-pullyu cases described above. 
Despite very different approaches to the proof of convergence to the traveling wave in the pushed and pulled cases, one 
can see one common feature in the original KPP results (\ref{e.c021601})-(\ref{mar2306})  
and in the pushed case (\ref{mar2308})-(\ref{mar2312}). Namely, the obtained rate of convergence of~$u(t,x)$ to $U_*(x)$ is 
much finer than the corresponding obtained rate of convergence for the front location. To see this, one needs to only compare  (\ref{e.c021601})
to (\ref{mar2306}) in the pulled case and~(\ref{mar2308}) to (\ref{mar2312}) in the pushed case.

Here,
we recover and explain this philosophy that ``rough front location asymptotics gives a finer
rate of convergence to a traveling wave."
We introduce a novel approach to quantifying the convergence rate in~\eqref{e.c021601} that provides one
simple explanation both for the exponential and algebraic rates in the pushed and pulled cases, respectively.  
Roughly, we prove the following (cf.~\Cref{t.rate-hr}), under some technical assumptions:
\be\label{e.rate}
	|u(t,m(t) + \cdot) - U_*(\cdot)| = \begin{cases}
				O(t^{-1}) \qquad&\text{ if } c_* = 2,\\
				O\big(\exp\big(-\farc{(c_*^2-4)t}{4}\big)\big)
						\qquad&\text{ if } c_* > 2.
			\end{cases}
\ee
As we have mentioned, in the case $c_*=2$, the convergence rate in (\ref{e.rate}) has been established in~\cite{Graham} for the 
Fisher-KPP nonlinearities based on the very fine asymptotics (\ref{mar2408}).  The proof here is completely different and avoids (\ref{mar2408})
altogether. For the other pulled and pushmi-pullyu cases the rate in~(\ref{e.rate}) is, to the best of our knowledge, new, as is the explicit rate 
in the pushed case.

To explain the approach to the proof of the convergence rates in (\ref{e.rate}), we need to recall the notion of the shape defect function introduced in~\cite{AHR2}.  
It is well known that the traveling wave solutions to (\ref{e.rde}) are monotonically decreasing.  Thus, there is a $C^1(0,1)$  function $\eta(u)$ so that
\be\label{e.eta}
	- U_*' = \eta(U_*).
\ee
It is easy to see that 
\be\label{mar2618}
	\eta(u)>0\hbox{ for all $u\in(0,1)$}
		\quad\text{ and }\quad
		 \eta(0)=\eta(1)=0.
\ee
We call $\eta(u)$ the ``traveling wave profile function.'' We define the shape defect function to be
\be\label{e.sdf_def}
w(t,x)= - u_x(t,x) - \eta(u(t,x)).
\ee
This, in a sense, represents how close the solution $u(t,x)$ is to solving~\eqref{e.eta} and is a measure of the ``distance in shape'' between $u(t,x)$ and the profile  $U_*(x)$.  
A major advantage here is that we do not {\em a priori} need to know {\em which} shift of $U_*$ is the closest one in order to use $w$ to obtain bounds on~$u(t,x) - U_*(x)$.  Imprecisely, one finds that
\be\label{e.w_eps}
	w = O(\eps)
		\quad\text{ if and only if }\quad
	u = U_* + O(\eps)
\ee
where the second inequality holds up to the appropriate shift.  We note that 
related quantities were used in~\cite{FifeMcLeod_pp,MatanoPolacik, Polacik, Uchiyama}; see~\cite{AHR2} for a more detailed discussion.

The main idea of this work is to estimate~$w(t,x)$ directly through its evolution equation
\be\label{e.sdf}
w_t - w_{xx} = w(Q(u) + \eta''(u) w), 
\ee
where, by \cite[equation (4.1)]{AHR2},
\be\label{e.Q}
	Q(u) = \eta'(u) (c_* - \eta'(u)) + \eta(u) \eta''(u)
		\qquad\text{ for all } u \in (0,1),
\ee
and use that information to read off the rate of convergence of $u(t,x)$ to the traveling wave profile~$U_*(x)$.  
As we see below, the nonlinearity $Q(u)$ satisfies
\be\label{jun2302}
Q(0)=f'(0) = 1
\ee
and, for a large class of nonlinearities, we also have 
\be\label{e.c062301}
	Q(u)	\leq 1
		\qquad\text{ for all } u\in[0,1],
\ee
see \Cref{l.general_eta}. 
 
A key informal observation is that if $u(t,x)$ is a solution to (\ref{e.rde}), 
there is a ``phantom front'' location $m_w(t)$ that is far behind the true front $m(t)$ and is where the shape defect function~$w(t,x)$ 
``wants'' to have its front.  The phantom front location of $w$ can be read off its equation~\eqref{e.sdf}. Surprisingly,
the evolution of $w(t,x)$ in (\ref{e.sdf}) turns out to be ``Fisher-KPP-like," regardless of whether the solution $u(t,x)$ to
(\ref{e.rde})  itself is of the pushed, pulled or pushmi-pullyu nature.
This is the main and, to us, unexpected
unifying element of all three cases. The simple reason behind this pulled nature of $w(t,x)$ is that, because of (\ref{jun2302})-\eqref{e.c062301},
ahead of the front it satisfies 
\be\label{jun2304}
w_t\leq w_{xx}+w,
\ee
which is exactly the same linearized problem as for the Fisher-KPP equation.

The second new key  point is that the distance
\be\label{jun2012}
D(t)=m(t)-m_w(t)
\ee
between the true and the phantom fronts 
controls the rate of convergence in (\ref{e.rate}), once again, regardless of whether the front is pushed or pulled. 
More precisely,  at an informal level, the main result of this paper is that the convergence rate in (\ref{e.rate}) comes from the estimate 
\be\label{mar2602}
	|u(t,m(t) + \cdot) - U_*(\cdot)|
		\sim |w(t, m(t) + \cdot)|
		= |w(t, D(t) + m_w(t) + \cdot)|
		\sim \exp\Big(-D(t)-\farc{D^2(t)}{4t}\Big),
\ee
where the first approximation follows from~\eqref{e.w_eps} and the second comes from the ``Fisher-KPP like" nature of (\ref{jun2304}); see also~\eqref{mar2634}, below.
In particular, this explains why one needs only ``rough'' asymptotics for $m(t)$ and $m_w(t)$ to get an ``exponentially finer''
convergence rate in~(\ref{e.rate}).  
In order to pass from (\ref{mar2602}) to (\ref{e.rate}), we show that, as long as $f(u)$ satisfies (\ref{e.f})-(\ref{e.normalization}) and some additional technical assumptions, 
the front location and the phantom front location have the 
following behavior as $t\to+\infty$:
\be\label{mar2604}
\begin{alignedat}{3}
&m(t)=c_*t+O(1),
	~\quad &m_w(t)=2t-\farc{3}{2}\log t+O(1),
	~\quad &\hbox{ in the pushed case},\\
&m(t)=2t-\farc{1}{2}\log t+O(1),
	~\quad &m_w(t)=2t-\farc{3}{2}\log t+O(1),
	~\quad &\hbox{ in the pushmi-pullyu case,}\\
&m(t)=2t-\farc{3}{2}\log t+O(1),
	~\quad &m_w(t)=2t-\farc{5}{2}\log t+O(1),
	~\quad & \hbox{ in the pulled case}.
\end{alignedat}
\ee
Using~(\ref{mar2602}) 
and~(\ref{mar2604}) leads directly to (\ref{e.rate}). 

The asymptotics for $m(t)$ in (\ref{mar2604}) in all three cases is already known and to a better precision
than stated in (\ref{mar2604}),
with the pushmi-pullyu case analyzed recently in~\cite{AHR2} and formally predicted in~\cite{berestycki2017exact,Ebert-vanSaarlos,Leach-Needham}. 
Our main goal here is to explain what the phantom front location $m_w(t)$
is, how~(\ref{mar2602}) comes about, and how the asymptotics of $m_w(t)$ in (\ref{mar2604})
can be computed.  We emphasize that, unlike~\cite{Graham,NRR2} that analyzed the Fisher-KPP
case, 
we only use the~$O(1)$-precise asymptotics for $m(t)$
and not anything finer to get the convergence rates in~(\ref{e.rate}). 

In all of the three cases in (\ref{mar2604}), the analysis of the phantom front location~$m_w(t)$ for the shape defect function is based on typical techniques for the Fisher-KPP equations (pulled fronts).  This leads to the surprising conclusion that, for a large class of nonlinearities, 
the convergence of the shifted solution~$u(t,x+m(t))$ 
to~$U_*(x)$ is  
a pulled phenomenon, regardless of the pushed, pulled, or pushmi-pullyu character of the spreading of~$u(t,x)$ itself. 
The reader may notice 
that the phantom front asymptotics~$m_w(t)$ in~(\ref{mar2604}) has the Bramson form (\ref{mar2404}), which is a signature of the  pulled fronts, precisely
when $m(t)$ is {\em not} pulled.  On the other hand, in the pulled case it is the front
asymptotics~$m(t)$ itself that has  the Bramson asymptotics
(\ref{mar2404}), while  
the phantom front position~$m_w(t)$ has an extra $\log t$ delay
relative to this location. 
This will be explained below. 
Of course,
without such a delay between~$m(t)$ and~$m_w(t)$, we would have~$D(t)=O(1)$ and~(\ref{mar2602}) 
would be useless!

We hope to  
convince the reader that the scheme outlined above is exceedingly simple to put into practice, beyond the situations we consider in 
the present paper.  Once one starts to work directly with
the shape defect function $w(t,x)$ and has the intuition (\ref{mar2602}),  the convergence proof is straightforward. 
In particular, the sometimes heavy computations, such as in the proof of Lemma~\ref{l.w_pulled} below, should not obfuscate this basic fact.  
We do not consider more general problems here because our interest is in the simplest possible presentation to illustrate the meaning behind the convergence rates.

\subsubsection*{Organization of the paper}

To better illustrate the method,
we first focus on the
the ``Hadeler-Rothe'' family of nonlinearities~$f$ given by~\eqref{e.HR} below.
In \Cref{s.results},  we give a statement of our main result, \Cref{t.rate-hr}, which establishes~\eqref{e.rate} in this context. This section also contains an expanded discussion both of the proof and of the sharpness of our bounds.  The proof of \Cref{t.rate-hr}, given in \Cref{s.sdf_proposition}, relies on estimates of the shape defect function in \Cref{t.sdf-hr}, which are proved in \Cref{s.sdf_proof}.

In order to analyze the evolution equation~\eqref{e.sdf} for $w$, we require 
some properties of the traveling wave profile function $\eta(u)$
and the nonlinearity  $Q(u)$ that appears in (\ref{e.sdf}). 
They are established in \Cref{s.Q} in some generality, not just for the 
Hadeler-Rothe nonlinearities.  Following this, \Cref{s.general} contains 
an extension of the 
convergence rates~\eqref{e.rate} to the general case.  The key observation is that 
the 
proof of \Cref{t.rate-hr}
uses the particular form of the Hadeler-Rothe nonlinearities essentially only 
through these properties of $Q$ and $\eta$.  
General versions of \Cref{t.rate-hr} are formulated there,  
in \Cref{t.rate-general,t.rate-general2}.

\subsubsection*{Acknowledgements}

CH was supported by NSF grants DMS-2003110 and DMS-2204615. LR was supported by NSF grants DMS-1910023 and DMS-2205497 and by ONR grant N00014-22-1-2174.  JA and CH acknowledge support of the Institut Henri Poincar\'e (UAR 839 CNRS-Sorbonne Université), and LabEx CARMIN (ANR-10-LABX-59-01).

\section{Convergence rates for the Hadeler-Rothe nonlinearities}\label{s.results}

To fix the ideas in a simple setting, we look in detail
at the special class of the so-called 
Hadeler-Rothe nonlinearities. They have the form
\be\label{e.HR}
f(u)=(u-u^n)(1+\chi nu^{n-1}),
\ee
with some $n\ge 2$ and $\chi\ge 0$.  The traveling waves for such nonlinearities 
were discussed in detail in~\cite{HadelerRothe,murray2007mathematical} for $n=2$ and in~\cite{Ebert-vanSaarlos} for $n>2$.
The classical Fisher-KPP nonlinearity~$f(u)=u-u^2$ is a special case of (\ref{e.HR}) with $\chi=0$ and $n=2$. 

It was shown in~\cite{Ebert-vanSaarlos,HadelerRothe,murray2007mathematical} for nonlinearities of the form (\ref{e.HR}) 
that there is a pushed-to-pulled transition at $\chi=1$: 
\be\label{mar2622}
	c_*(\chi)
		= \begin{cases}
			2
				\qquad&\text{ if } 0 \leq \chi \leq 1,\\
			\sqrt \chi + \frac{1}{\sqrt \chi}
				\qquad&\text{ if } \chi \geq 1.
		\end{cases}
\ee
Moreover, the traveling wave profile function is explicit for $\chi\ge 1$ and is given by
\be\label{mar2621}
\eta(u) = \sqrt \chi (u-u^n),
\ee
see \cite[Proposition~A.2]{AHR2}.  Hence, when $\chi\ge 1$,
the traveling waves  have the purely exponential asymptotics~(cf.~\eqref{e.c062203}): there exists $\eps, A_1 > 0$ so that
\be\label{mar2610}
	U_*(x)\sim A_1 e^{-\lambda_0 x} + O(e^{-(\lambda_0 + \eps) x}),~~\hbox{as $x\to+\infty$}.
 \ee 
When $0 \leq \chi < 1$, no such
explicit expression is possible for $\eta(u)$ because $U_*$ has the pulled asymptotics:
there exists some $\eps>0$ and $A_0 >0$ so that
\be\label{e.lambda}
	U_*(x) \sim (A_0 x + B_0) e^{- \lambda_0 x} + O( e^{-(\lambda_0 +\eps) x}),~~\hbox{as $x\to+\infty$}.
 \ee
The decay rate $\lambda_0>0$ in~\eqref{mar2610} and~\eqref{e.lambda} is the largest root of
\be\label{mar2814}
	c_* \lambda_0 = \lambda_0^2 + f'(0).
\ee
Recalling~\eqref{e.normalization}, if $c_* = 2$, then $\lambda_0 = 1$.
Let us mention that, after a spatial shift, we may assume that $B_0=0$,
so that (\ref{e.lambda}) becomes
\be\label{july1002}
U_*(x) \sim A_0xe^{- \lambda_0 x} + O( e^{-(\lambda_0 +\eps) x}),~~\hbox{as $x\to+\infty$}.
 \ee
This is another natural normalization that we will sometimes 
use below as an alternative to~(\ref{jun2004}).

The corresponding front location asymptotics for the solutions to~(\ref{e.rde})
with a rapidly decaying initial condition was established in~\cite{AHR2}: there exists $x_0$ that depends on the initial condition $u_0$,
so that, as $t\to\infty$
\be\label{mar2710}
\begin{alignedat}{2}
m(t)&=2t-\farc{3}{2}\log t+x_0,  \quad\qquad&&\text{for $0\le\chi<1$ (the pulled case),}\\
m(t)&=2t-\farc{1}{2}\log t+x_0,  \quad\qquad&&\text{for $\chi=1$ (the pushmi-pullyu case),}\\
m(t)&=c_*(\chi)t+x_0,  \quad\qquad &&\text{for $1<\chi$ (the pushed case).}
\end{alignedat}
\ee

It is convenient to recall the asymptotic behavior of $U_*$ as $x \to -\infty$ as well: there are $A_1, \eps > 0$ so that
\be\label{mar2614}
1 - U_*(x)\sim
A_1 e^{\lambda_1 x} +  O(e^{(\lambda_1+\eps) x}),~~\hbox{ as $x \to -\infty$.}
\ee
Here, $\lambda_1$ is the nonnegative root of
\be\label{e.lambda_bis}
- c_* \lambda_1 = \lambda_1^2 + f'(1).
\ee
Notice that, due to~\eqref{e.HR}, we have
\be\label{e.c062401}
	\lambda_1 > 0
		\qquad \text{ since }\quad
	f'(1) = -(n-1) (1 + \chi n) < 0.
\ee

\subsection{The main result for the Hadeler-Rothe nonlinearities} 

In this section, we state the convergence rates in (\ref{e.rate}) for the Hadeler-Rothe nonlinearities of the form (\ref{e.HR}).
For simplicity, we take an initial condition $u(0,x) = u_0(x)$ such that~$0\le u_0(x)\le 1$ for all $x\in\Rm$, and there exsts some $L_0\in\R$,
so that
\be\label{e.u_0}
	u_0(x) = 0
		 \text{ if } x \geq L_0,~~
	 \text{and} ~~
	w_0 (x)= w(0,x) \geq 0,~~\hbox{ for all $x\in\Rm$.}
\ee
The non-negativity assumption on $w(0,x)$ simply encodes that the initial condition $u_0(x)$ is ``steeper'' than $U_*(x)$.  
In particular, it follows from (\ref{e.u_0}) that $u_0(x)$ is monotonically decreasing and  $u_0(x)\to 1$ as~$x\to-\infty$. 
The comparison principle and~\eqref{e.sdf} yield
that then $u(t,x)$ remains steeper than $U_*(x)$ for all $t>0$,  in the sense that
\be\label{e.w>0}
w(t,x) > 0,~~\text{ for all } t> 0,~x \in \R.
\ee
A typical example of such initial condition is $u_0(x)=\1(x\le 0)$.  
We believe that the non-negativity assumption on $w(0,x)$ can be relaxed by using results such as by Angenent in~\cite{Angenent} 
or Roquejoffre in~\cite{Roq-monot} to show that $w(t,x)$ 
``eventually'' becomes nonnegative, at least on every compact set.  
We adopt this assumption to avoid the related technicalities.

Our main result for the Hadeler-Rothe nonlinearities is as follows.
\begin{thm}\label{t.rate-hr}
Suppose that $u$ solves~\eqref{e.rde} with nonnegative initial condition $u_0$ satisfying~\eqref{e.u_0}.  
Assume that $f(u)$ is given by (\ref{e.HR}) with some $\chi\ge 0$ and $n\ge 2$.  Let $c_*$ be given by (\ref{mar2622}). 
Then there is $\sigma: [0,\infty) \to \R$
so that:
\begin{enumerate}[(i),itemindent=-7pt]
	\item if~$0\le\chi\le 1$, then
\be\label{mar2806}
\|u(t,\cdot + \sigma(t)) - U_*(\cdot)\|_{L^\infty} \leq \frac{C}{t},
\ee

	\item if $\chi>1$, then for any $\Lambda >0$,
		\be\label{e.c81901-hr}
			\|u(t,\cdot + \sigma(t)) - U_*(\cdot)\|_{L^\infty([-\Lambda,\infty))}
				\leq \frac{C_\Lambda}{\sqrt t} e^{- \frac{(c_*^2- 4)}{4} t }.
		\ee
\end{enumerate}
\end{thm}
As will be seen from the proof, convergence occurs in a (stronger) weighted $L^\infty$-norm, 
but we opt for the simpler statement here.

The main ingredients in \Cref{t.rate-hr} are knowledge of the true front location $m(t)$ as well as the behavior of $Q$ and $\eta$ in \eqref{e.c062301}.  In this sense, we use the form~\eqref{e.HR} in a rather weak way.  We provide a full discussion of the general case in \Cref{s.general} and formulate broader versions of \Cref{t.rate-hr} there; see \Cref{t.rate-general,t.rate-general2}.

Interestingly, unlike the classical results in~\cite{fife1977approach,FifeMcLeod_pp,Rothe} for pushed
waves, the estimate (\ref{e.c81901-hr}) does not depend on $f'(1)$.  Actually, a similar argument using our methods yields a messier global estimate:
\be\label{e.c81902}
\|u(t,x + \sigma(t)) - U_*(x)\|_{L^\infty}
\leq Ce^{- \min\big(\frac{c_*^2- 4}{4},|f'(1)|\big) t + o(t)}.
\ee
However, the $f'(1)$ term in the exponential merely reflects the ``slowness'' with which $U_*$ converges to~$1$ on the left.  
We choose to present the ``at and beyond the front'' estimate (\ref{e.c81901-hr}) 
above because it is a better representation of the mechanism that pulls $u(t,x)$ towards $U_*(x)$.  
In particular, it reflects the aforementioned pulled nature of the convergence of the solution to the wave in shape, regardless of whether the wave itself is pushed or pulled.

\subsection{Discussion of the proof}\label{s.discussion}

A very useful observation is that, for the Hadeler-Rothe nonlinearities,~\eqref{e.c062301} holds and the traveling wave profile function $\eta(u)$ is concave.   
\begin{prop}\label{prop-mar2602}
Assume that $f(u)$ has the form (\ref{e.HR}), then, for any $\chi \geq 0$ and $n\geq 2$,
\be
	Q(u) \leq 1
		\quad\text{ and }\quad
	\eta'' (u)\leq 0, ~~\text{ for all $u \in (0,1)$.}
\ee
\end{prop}
A more precise version is stated in \Cref{l.R_bound}.
\Cref{prop-mar2602} follows immediately from the explicit expression (\ref{mar2621}) for $\eta(u)$ when $\chi \geq 1$.  Otherwise, it is proved in \Cref{l.R_bound}.  Its generality, beyond the Hadeler-Rothe class, is discussed in \Cref{s.general}.

\Cref{prop-mar2602} is nearly enough to understand the phantom front $m_w(t)$ as we have, at highest order,
\be\label{jul1004}
	w_t \approx w_{xx} + w
\ee
ahead of the front.  Remarkably, this is exactly the same as the linearization for the classical Fisher-KPP equation
\be
u_t=u_{xx}+u-u^2.
\ee 
This would suggest that $m_w(t)$ should be given 
by the standard Bramson asymptotics~\eqref{mar2404} for the Fisher-KPP case. 
However, it has been observed that 
the Bramson shift may be sensitive to lower order terms ahead of the 
front for nonlinearities that are not better than Lipschitz near 
$u=0$~\cite{Bou-Hen}. In that case, (\ref{jul1004}) may be not a faithful approximation
to (\ref{e.sdf}).  
 It is, thus, crucial to understand the regularity of $\eta$ near $u = 0$.  As a consequence, we consider two cases depending on this regularity.

\subsubsection*{The pushed and pushmi-pullyu cases: $\chi \geq 1$}
Consider first the pushed and pushmi-pullyu cases, where $\eta$ is given 
explicitly by~\eqref{mar2621} and is smooth at $u=0$.  In this case,
\be\label{mar2682}
\bal
Q(u)&
=1-n(1-2\chi+\chi n )u^{n-1}-\chi n u^{2n-2}
= 1 + O(u^{n-1})
	\qquad \text{ as $u\to 0$}.
\enbal
\ee
Recall that $n\geq 2$.  Hence, we expect that, ahead of the front of $u(t,x)$, 
the shape defect function~$w(t,x)$ does behave 
approximately as a solution to
\be\label{mar2631}
w_t=w_{xx}+w,
\ee
when $\chi\ge 1$.  An informal consequence of~\cite{HNRR} 
is that $w(t,x)$, being bounded and approximately satisfying~\eqref{mar2631} where it is small, 
``wants to have a front'' at the location 
\be\label{mar2633}
m_w(t)=2t - \frac32 \log t,
\ee
and should have the approximate form
\be\label{mar2634}
w(t, x + m_w(t))
		\approx \exp\Big\{ - x - \frac{x^2}{4t} + \text {(lower order terms)}\Big\},~~\hbox{ for $x\gg 1$.}
\ee
On the other hand, $w(t,x)$ is governed by $u(t,x)$, which has its front at the position $m(t) = c_*t$ in the pushed case~$\chi>1$, 
and at $m(t) = 2t - \sfrac12 \log t$ in the pushmi-pullyu case
$\chi=1$~\cite{AHR2}.  
Hence, we have, up to lower order terms
\be
	D(t)
		= m(t) - m_w(t)
		\approx \begin{cases}
			\log t
				\qquad&\text{ if } \chi = 1,\\
			(c_* - 2)t
				\qquad&\text{ if } \chi > 1.
		\end{cases},
\ee
According to (\ref{mar2634}), this produces
\be\label{mar2635}
	w(t,m(t))
		= w(t, D(t) + m_w(t))
		\approx \exp\Big\{-D(t) - \frac{D^2(t)}{4t}\Big\},
\ee
which, along with~\eqref{e.w_eps}, yields \Cref{t.rate-hr}.

Let us note that the explicit form of $\eta$, beyond \Cref{prop-mar2602}, is not needed here, because the key estimate used above, that is, the right hand side of~\eqref{mar2682}, follows directly from the traveling wave asymptotics~\eqref{mar2610} and~\eqref{e.f_eta} below.  Indeed, we can see that, whenever~\eqref{mar2610} holds, we have, for some $\alpha>0$,
\be\label{e.eta_asymp_pushed}
	\eta(u) \sim u + O(u^{1+\alpha}).
\ee
See \Cref{lem.eta_pushed}.

\subsubsection*{The pulled case: $0\leq \chi < 1$}

For $0\le\chi<1$, we do not have an explicit expression for $\eta(u)$ or $Q(u)$.
To understand the behavior of $Q(u)$ for $u\ll 1$ in this range of $\chi$, we can, at least informally, deduce the behavior of $\eta$ and its derivatives from~\eqref{e.lambda}.

Using~\eqref{e.eta}, we can write two useful identities involving $\eta$:
\be\label{e.f_eta}
	f(u) = \eta(u)(c_* - \eta'(u))
	\qquad\text{ and }\qquad
	\eta(u) = - U_*'\circ U_*^{-1}(u).
\ee
From these, we immediately observe that
\be\label{e.eta_C2}
	\eta \in C^\infty_{\rm loc}(0,1),
	\qquad
	\eta'(0) = \lambda_0,
	\qquad\hbox{ and }\qquad
	\eta'(1) = -\lambda_1.
\ee
Both~\eqref{e.f_eta} and~\eqref{e.eta_C2} hold for any $f$ satisfying~\eqref{e.f}-\eqref{e.normalization}.  The endpoint regularity is more subtle and is affected by the  additional linear factor in~\eqref{e.lambda} that is present in the pulled case.
Indeed, from~\eqref{e.lambda}, it is straightforward to see that
\be\label{e.eta_asymp_pulled}
	\eta(u)\sim u+\farc{u}{\log u},\qquad \hbox{as $u\to 0$},
\ee
from which we formally deduce that
\be\label{e.eta'}
	\eta'(u)\sim 1+\farc{1}{\log u}
		\qquad\text{ and }\qquad
	\eta''(u)\sim-\farc{1}{u\log^2u}
	\qquad\text{ as } u \to 0^+.
\ee
These are made 
precise in \Cref{l.eta_pulled} below. 
Therefore, when $0\le\chi<1$, the function $Q(u)$ defined in (\ref{e.Q}) has the asymptotics
\be\label{mar2630}
Q(u)\sim 1-\farc{2}{\log^2u},~~\hbox{as $u\to 0$}.
\ee 
Thus, 
a good approximation to $w(t,x)$ is by a solution to a modification of (\ref{mar2631}):
\be\label{e.c051701}
w_t - w_{xx} \approx w\Big(1 - \frac{2}{\log^2u}\Big).
\ee
Using, once again very informally, the main result of~\cite{Bou-Hen}, we see that 
the shape defect function~$w(t,x)$ ``wants to have its front'' at the location
\be\label{mar2636}
m_w(t)=2t - \frac52 \log t,
\ee
while the front of $u(t,x)$ is at the Bramson position 
\be
m(t)=2t - \frac32 \log t,
\ee
as follows from~\cite{AHR2}.  Thus, for $0\le\chi<1$, we have
$D(t) = \log t$ and~\eqref{mar2635} again yields the $O(\sfrac1t)$ convergence 
rate in~\eqref{e.rate}.

The above informal arguments indicate that, as we have already mentioned, 
the behavior of the shape defect function $w(t,x)$ is 
always a pulled phenomenon regardless of the pushed, pulled, or pushmi-pullyu spreading of $u(t,x)$ itself.

\subsection{Sharpness of \Cref{t.rate-hr}}

It appears that this approach leads to matching lower bounds.  This is easiest to see in the pushed case.  Indeed, fixing $\eps, \delta \ll 1$, $R\gg 1$, and $T\gg 1$, it is 
straightforward  to check that
\be
	\underline w(t,x+ (c_* + \eps) t)
		= \delta e^{ t\left(1 - \frac{\pi^2}{4R} - \frac{c_*^2}{4} - C\delta\right)} e^{- \frac{c_* + \eps}{2} x} \cos\Big(\frac{x \pi }{2 R}\Big) \1_{[-R,R]}(x)
\ee
is a subsolution to~\eqref{e.sdf} for $t\geq T$. The additional $\eps t$ shift in 
the moving frame allows us to use the approximation~$Q \approx 1$ 
because it puts us in the regime where $u \ll 1$.  Up to further adjusting $\delta$, it is easy to check that~$\underline w(1, \cdot) \leq w(1,\cdot)$.  It follows that
\be
	e^{ t\left(1 - \frac{\pi^2}{4R} - \frac{c_*^2}{4} - C\delta\right)}
		\leq C w(t, x + (c_* + \eps) t)
			\qquad\text{ for all } x\in [-\sfrac{R}{2},\sfrac{R}{2}].
\ee
From this, a simple ODE argument shows that
\be
	\|u(t,\cdot + \sigma(t)) - U_*(\cdot)\|_{L^\infty} \geq e^{- \frac{(c_*^2-4)}{4} t + o(t)}.
\ee

The arguments in the pulled and pushmi-pullyu cases will be more involved.  We, nonetheless, expect them to proceed in a fairly straightforward manner using the shape defect function.

\section{Estimates on the shape defect function}\label{s.sdf_proposition}

One of the main technical points of this paper is that 
the proof of Theorem~\ref{t.rate-hr} requires understanding the front location asymptotics for $u(t,x)$ only up to $O(1)$ as $t\to+\infty$. 
For the Hadeler-Rothe nonlinearities we have the following.  
\begin{prop}[\cite{AHR2}]\label{p.weak_front}
Under the assumptions of \Cref{t.rate-hr}, let the function $m(t)$ be given by~(\ref{mar2710}). Then, we have
\be\label{mar2802}
\lim_{L\to\infty} \limsup_{t\to\infty} \sup_{x \geq m(t)+L} u(t,x) = 0
\quad\text{ and }\quad
\lim_{L\to\infty} \liminf_{t\to\infty} \inf_{x \leq m(t)-L} u(t,x) = 1.
\ee
\end{prop}
This claim holds, of course, for a much wider class of nonlinearities -- see~\cite{AHR2,Giletti} for a discussion. 
The next lemma gives preliminary control on how quickly $u(t,x)$ 
tends to its limits as~$x\to\pm \infty$.      
\begin{lem}\label{l.steepness}
With $m(t)$ as in \Cref{p.weak_front} and $w(t,x)$ satisfying~\eqref{e.w>0}, there is $C>0$ so that
\be
	u(t,x + m(t))
		\geq U_*(x + C) \quad\text{ for all } x < 0,
	\quad\text{ and }\quad
	u(t,x + m(t))
		\leq U_*(x - C) \quad\text{ for all } x > 0.
\ee
\end{lem}
By a simple ODE comparison argument
using~\eqref{e.eta},~\eqref{e.sdf_def}, and~\eqref{e.w>0}, we see that, for any $x_1,x_2$, 
\be\label{e.c062701}
	\text{ if } u(t,x_1) = U_*(x_2)
		\quad\text{ then }\quad
		u(t, x_1 + x)
			 \begin{cases} 
				\leq U(x_2 + x) &\quad\text{ if } x > 0,\\
				\geq U(x_2 + x) &\quad\text{ if } x < 0.
			\end{cases}
\ee
Then \Cref{l.steepness} follows directly from \Cref{p.weak_front}.  The proof is omitted.

The main step allowing us to deduce the bounds in Theorem~\ref{t.rate-hr} 
is the following estimate on the shape defect function at the front location $m(t)$.  
\begin{thm}\label{t.sdf-hr}
Suppose the assumptions of \Cref{t.rate-hr} hold.  Let $m(t)$ and $\lambda_1>0$ be as in~(\ref{mar2710}) and~\eqref{e.lambda_bis}, respectively, and let $\eps >0$.\\
(i) If $0\le\chi<1$, then %
 \be
w(t,x + m(t))
\leq \frac{C}{t} \Big((1 + x)^2 e^{-x - \frac{x^2}{Ct}}\Big)\one(x\ge 0)
+ \frac{C_\eps}{t} 
	e^{(\lambda_1-\eps) x}\one(x\le 0).
\ee
(ii) If $\chi=1$ then 	 		
\be
w(t,x + m(t)) \leq \frac{C}{t} \Big((1 + x) e^{- x - \frac{x^2}{Ct}}\Big) \one(x\ge 0)
	+ \frac{C_\eps}{t}
		e^{(\lambda_1-\eps) x}\one(x\le 0).
\ee
(iii)  If $\chi>1$ and $x > L_0-m(t)$ (recall $L_0$ from~\eqref{e.u_0}) then 
\be\label{e.w.pushed-hr}
	w(t,x + m(t)) \leq \frac{C}{\sqrt t} 
	\exp\left\{ - \frac{c_*^2 -4}{4} t - \frac{c_*x}{2} - \frac{x^2}{4t}\right\},
\ee
with $c_* = c_*(\chi)$ given by (\ref{mar2622}).
\end{thm}

We note that the $\eps$ in cases (i) and (ii) can almost certainly be removed with a more careful proof.  Our focus in this paper, 
however, is not on the sharpest possible behavior on the left, as~$x\to-\infty$.  

While the statements in Theorem~\ref{t.sdf-hr}(i)-(ii) for the  pulled 
and pushmi-pullyu cases are slightly 
different, the proofs, postponed until \Cref{s.sdf_proof}, are nearly identical.  They are based on the intuition discussed in \Cref{s.discussion}: 
the equation for~$w(t,x)$ wants to spread slower than the equation for $u(t,x)$. 
The statement of \Cref{t.sdf-hr}(iii) in the pushed case and its proof, presented in Section~\ref{subsec:pushed-proof},
are different because we can use an elementary estimate ``out-of-the-box''.

\subsection{Deducing \Cref{t.rate-hr} from \Cref{t.sdf-hr}} \label{s.main_thm}

\subsubsection{Preliminary bounds on $\eta$}

We now make the behavior of $\eta(u)$ near $u=0$, stated informally in~\eqref{e.eta'}, precise.   
\begin{lem}[Asymptotics of $\eta(u)$ in the pulled case]\label{l.eta_pulled}
	Assume that $f\in C^2([0,1])$ and satisfies~(\ref{e.f})-(\ref{e.normalization}).  Suppose that the profile $U_*(x)$  has the asymptotics (\ref{e.lambda}) 
	as $x\to+\infty$.  Then there exists $C>0$ so that, for $u\in (0, \sfrac1{100})$,
\begin{multicols}{2}
	\begin{enumerate}[(i)]
		\item \quad$\displaystyle\Big|\eta(u) - \Big(u + \frac{u}{\log u}\Big)\Big|
				\leq C\frac{u \log \log(\sfrac1u)}{\log^2(\sfrac1u)},$
		\item \quad$\displaystyle\Big|\eta'(u) - \Big(1 + \frac{1}{\log u}\Big)\Big|
				\leq C\frac{\log \log(\sfrac1u)}{\log^2(\sfrac1u)},$
				
		\item \quad$\displaystyle\Big|\eta(u)\eta''(u) - \Big(\frac{-1}{\log^2 u}\Big)\Big|
				\leq C\frac{\log \log(\sfrac1u)}{\log^3(\sfrac1u)}.$
				
		\item[]
	\end{enumerate}
\end{multicols}
\end{lem}	
We note that this lemma does not require 
the specific form~\eqref{e.HR} of $f$.  The parts~(i)-(ii)  
will be used to deduce \Cref{t.rate-hr} from \Cref{t.sdf-hr}.  The property (iii)  is not required for that proof but will be needed in the
proof of   \Cref{t.sdf-hr} itself.

\begin{proof}
We use the normalization of $U_*(x)$ in which $B_0 = 0$ in~\eqref{e.lambda}. 
Consider first the claim (i). Fix $u\in(0,\sfrac1{100})$ and $x_u$ such that $U_*(x_u) = u$.  
We deduce from (\ref{e.lambda}) with $B_0=0$ that  
\be\label{e.c051602}
	x_u = \log\sfrac1u + O\Big(\log\log\sfrac1u\Big),~~\hbox{as $u\to 0^+$.} 
\ee
Using this in the definition of $\eta(u)$, we find
\be\label{e.c051603}
	\begin{split}
	\eta(u)
		&= \eta(U_*(x_u))
		= - U_*'(x_u)
		= A_0 x_u e^{-x_u} - A_0 e^{-x_u} + O(e^{-(1+\eps)x_u})
		\\&
		= U_*(x_u) \left(1 - \frac1{x_u} + O(x_u^{-1}e^{-\eps x_u})\right),
	\end{split}
\ee
The claim (i) follows then from inserting~\eqref{e.c051602} into~\eqref{e.c051603} and using a straightforward expansion.

We omit the proofs of (ii) and (iii) as they proceed by similar arguments.
\end{proof}

\begin{lem}[Asymptotics of $\eta$ in the pushed and pushmi-pullyu cases]\label{lem.eta_pushed}
Assume that $f\in C^2([0,1])$ and satisfies~(\ref{e.f})-(\ref{e.normalization}).  Suppose that the profile $U_*$ has the asymptotics (\ref{mar2610}) as $x\to+\infty$.  
Then, there exist $\alpha>0$ and $C>0$ such that, for all $u\geq 0$,
\be\label{mar2808} 
	|\eta'(u) - \lambda_0| \leq Cu^\alpha.
\ee

\end{lem}

The proof is omitted as it is a simpler version of the proof of Lemma~\ref{l.eta_pulled}.

\subsubsection{The proof of \Cref{t.rate-hr}}

The first steps of the proof for both cases (i) and (ii) can be handled simultaneously. 
As $u(t,x)$ is monotonic in $x$, we may define~$\sigma(t)$   by
\be\label{mar2716}
u(t, \sigma(t))= U_*(0).
\ee
We shift to the corresponding moving frame: let
\be
	\tilde u(t,x) = u(t, x+ \sigma(t))
	\quad\text{ and }\quad
	\widetilde w(t,x)=u(t,x+\sigma(t)).
\ee
It follows from \Cref{p.weak_front} that
\be\label{e.c60106}
	\sup_{t \geq 1} |\sigma(t) - m(t)| \leq C.
\ee
We may then apply \Cref{t.sdf-hr} with $\sigma(t)$ in place of $m(t)$, at the expense of changing the constants. 

To use Theorem~\ref{t.sdf-hr}, we need to bound the smallness of the difference
\be
s(t,x) = \tilde u(t,x) - U_*(x)
\ee
in terms of the smallness of the shape defect function $\widetilde w(t,x)$. 
Note that, by the choice of $\sigma(t)$ in~(\ref{mar2716}), 
\be\label{mar2718}
s(t,0)=0,~~\hbox{ for all $t>0$.}
\ee
We also point out that by the steepness comparison~\eqref{e.c062701}, 
we have
\be\label{e.c81801}
 s(t,x) \leq 0 \text{ when } x >0,
	~~\text{and}~~
	s(t,x) \geq 0 \text{ when } x < 0.
\ee

In order to relate $s(t,x)$ to $\widetilde w(t,x)$, note that, for each fixed $t$, $s(t,x)$ satisfies the following ODE in $x$:
\be\label{e.c60101}
s_x
		= -\widetilde w - \eta(\tilde u) + \eta(U_*)
		= -\widetilde w - \eta'(\xi(t,x)) s.
\ee
Here, $\xi(t,x)$ is  an intermediate point between $\tilde u(t,x)$ and $U_*(x)$ given by the mean value theorem.  
From~\eqref{e.c60101}, we obtain 
\be\label{e.c60205}
	\Big(\exp\Big\{\int_0^x \eta'(\xi(t,y)) dy\Big\} s(t,x)\Big)_x
		= -\exp\Big\{\int_0^x \eta'(\xi(t,y)) dy\Big\} \widetilde w(t,x).
\ee
Using the boundary condition (\ref{mar2718}) and integrating gives 
\be\label{e.c60104}
	\begin{split}
	s(t,x)
	&= - \exp\Big\{-\int_0^x \eta'(\xi(t,z)) dz\Big\} \int_0^x \exp\Big\{\int_0^y \eta'(\xi(t,z)) dz\Big\} \widetilde w(t,y) dy
			\\&
			= - \int_0^x \exp\Big\{- \int_y^x\eta'(\xi(t,z)) dz\Big\} \widetilde w(t,y) dy.
	\end{split}
\ee

From here, the main points of the proof are exactly the same in each case (i)-(ii); however, due to the difference in the precise asymptotics in \Cref{t.sdf-hr} in these
two cases, we have no choice but to write up each case separately.

\subsubsection*{Proof of \Cref{t.rate-hr}(i)} 
We analyze separately the cases~$\chi\in[0,1)$ and $\chi=1$,
as in parts  (i) and (ii) of \Cref{t.sdf-hr}.

Fix $\chi \in [0,1)$.  We consider first  $x \geq 0$, so that  $s(t,x)\le 0$ due to~\eqref{e.c81801}.  Thus, we only have to obtain a lower bound on $s(t,x)$.    
In view of
(\ref{e.c60104}), 
we seek control on the $\eta'(\xi(t,x))$ term.   We have
\be\label{e.c81802}
\tilde u(t,x)\le \xi(t,x)\leq U_*(x)\leq C (x+1) e^{-x}.
\ee
Using~\eqref{e.c81802} and the asymptotics  in \Cref{l.eta_pulled}(ii), gives, if $x$ is sufficiently large,
\be\label{jun2010}
\bal
\eta'(\xi(t,x))
&	\geq 1 + \frac{1}{\log \xi} - \frac{C\log\log \frac{1}{\xi}}{\log^2\xi}\geq 1 + \frac{1}{(-x)+\log x+\log C} - \frac{C}{(x+1)^{\sfrac32}}\\
&\ge   1 + \frac{1}{(-x) } - \frac{C}{(x+1)^{\sfrac32}}.
\enbal
\ee
Hence, for all $x\geq 0$, we have, after increasing the constant $C>0$ in (\ref{jun2010}), 
\be\label{e.c60105}
	\eta'(\xi(t,x))
		\geq 1 - \frac{1}{x+ 1} - \frac{C}{(x+1)^{\sfrac32}}.
\ee  

Using~\eqref{e.c60105} and \Cref{t.sdf-hr}(i) in~\eqref{e.c60104}, with $m(t)$ replaced by $\sigma(t)$, yields
\be
	\begin{split}
	s(t,x)
		&\geq  - C \int_0^x \frac{x+1}{y+1} e^{-(x-y)} \tilde w(t,y) dy
		\geq - \frac{C (x+1) e^{-x}}{t} \int_0^x (y+1) e^{ - \frac{y^2}{Ct}} dy
		\\& 
		\geq - C (x+1)e^{-x} \min\bigg\{\frac{(x+1)^2}{t}, 1\bigg\}.
	\end{split}
\ee
This concludes the proof of (\ref{mar2806}) in the pulled case $0\le\chi<1$ for $x\ge 0$.

Next, consider the pushmi-pullyu case $\chi=1$, corresponding to \Cref{t.sdf-hr}(ii), still for $x\ge 0$.  Here, we can replace~\eqref{e.c81802} and~\eqref{e.c60105}    
by, respectively,  
\be
	\xi(t,x)\leq U_*(x)\leq C e^{-x}
		\qquad\text{ and }\qquad
	\eta'(\xi(t,x))\geq 1 - \frac{C}{x^{1+\alpha}},
\ee
due to the asymptotics in \Cref{lem.eta_pushed}.  Arguing similarly as above, we obtain
\be
	\begin{split}
	s(t,x)
		&\geq  - C \int_0^x e^{-(x-y)}\tilde  w(t,y) dy
		\geq - \frac{C e^{-x}}{t} \int_0^x (y+1) e^{ - \frac{y^2}{Ct}} dy
		\geq - Ce^{-x} \min \Big\{1, \frac{(x+1)^2}{t}\Big\}.
	\end{split}
\ee
This concludes the proof of \Cref{t.rate-hr}(i) on the domain $x \geq 0$ for $0\le\chi\le 1$.

We now consider the case $x\leq 0$.  Due to~\eqref{e.c81801}, we need only obtain an upper bound on $s(t,x)$.  The  argument is  essentially the same as for $x\ge 0$.  
The main differences are the asymptotics of $\eta'(u)$ near $u\approx 1$ and $U_*(x)$ and $w(t,x)$ as $x \to-\infty$.  Unlike before, we need not separate into the two cases, as the behavior at the back is the same both for $0\le\chi<1$ and $\chi=1$.

First, notice that
\be\label{e.c81803}
	1- C e^{\lambda_1 x} \leq U_*(x) \leq \xi(t,x)
		\qquad\text{ for } x \leq 0,
\ee
and, for all $u > 1/2$,
\be
	\eta'(u) = - \lambda_1 + O((1-u)^p)
		\quad\text{ for some } p > 1,
\ee
The combination of these two inequalities   leads to
\be\label{e.c81804}
	\eta'(\xi(t,x))
		\geq - \lambda_1 - C e^{- \eps x},~~\hbox{ for $x\le 0$},
\ee
where $\eps$ is as in~\eqref{mar2614}.

We use~\eqref{e.c60104} and then~\eqref{e.c81804} and  \Cref{t.sdf-hr}(i)-(ii) to find
\be
	\begin{split}
		s(t,x)
			&\leq\frac{C_\eps}{t} \int_x^0 e^{\lambda_1 (x-y)} e^{(\lambda_1-\eps) y} dy
			\leq \frac{C_\eps}{t} e^{\lambda_1 x} \int_x^0 e^{- \eps y} dy
			\\&
			= \frac{C_\eps}{t\eps} (e^{(\lambda_1 - \eps) x} -e^{\lambda_1 x})
			\leq \frac{C_\eps}{t} e^{(\lambda_1 - \eps) x} ,~~\hbox{ for $x\le 0$}.
	\end{split}
\ee
Here $\eps\in (0,\lambda_1)$ is arbitrary. 
This completes the proof of \Cref{t.rate-hr}(i).

\subsubsection*{Proof of \Cref{t.rate-hr}(ii)}
We proceed as above.  By the Harnack inequality, it suffices to consider the case $L=0$, so that~$x\ge 0$. 
Again, due to~\eqref{e.c81801}, we need only establish a lower bound on $s(t,x)$.  
Next, note that, due to \Cref{lem.eta_pushed}, we have, for some $p>1$,
\be
	\eta'(\xi(t,x))
		\geq \lambda_0 - \frac{C}{(1 + x )^p}.
\ee
We find, from~\eqref{e.c60104} and \Cref{t.sdf-hr}(iii), once again, with $m(t)=c_*t+x_0$, 
\be
	\begin{split}
	s(t,x)
		&\geq - C\int_0^x e^{-\lambda_0 (x-y)} \tilde w(t,y) dy
		\geq - \farc{C}{\sqrt t} e^{- \frac{c_*^2 - 4}{4} t} \int_0^x 
		e^{-\lambda_0(x-y) - \frac{yc_*}{2} - \frac{y^2}{4t}} dy
		\\&
		\geq  - \frac{C}{\sqrt t} e^{- \frac{c_*^2 - 4}{4} t} \int_0^x e^{-\lambda_0(x-y) - \frac{yc_*}{2}} dy
		\geq - \frac{C}{\sqrt t} e^{- \frac{c_*^2 - 4}{4} t} e^{- \frac{xc_*}{2}}.
	\end{split}
\ee
The second to last equality uses that $\exp\{- \sfrac{y^2}{4t}\} \leq 1$ and the last inequality  
uses that $\lambda_0 > c_*/2$, which follows from~\eqref{mar2814}. 
This concludes the proof. $\Box$

\section{The proof of \Cref{t.sdf-hr}}\label{s.sdf_proof}

Before we begin, we state one final lemma about the behavior of $\eta$ and $Q$, defined in~\eqref{e.eta} and~\eqref{e.Q}, respectively.  
This is the key and essentially only place 
in this manuscript where we use the form~\eqref{e.HR} of the Hadeler-Rothe nonlinearities $f(u)$. 
\begin{lem}\label{l.R_bound}
Suppose the assumptions of \Cref{t.rate-hr} hold. 
Then
\be\label{e.c051604}
	\eta'' (u)\leq 0
	\quad\text{ and }\quad
	Q (u)\leq 1,~~\hbox{ for all $u\in(0,1)$.} 
\ee
Further, we have the refined bounds: letting
\be\label{mar2902}
	\Rd(u)
		= 1- Q(u(t,x)), 
\ee
for any $\delta_0,\delta_1\in(0,\sfrac{1}{100})$ with $\delta_1$ sufficiently small, there are $r_0>0$ and $r_1 > 0$ such that
\be\label{mar2812}
\Rd(u) \ge \begin{cases}
			r_0,
				\qquad&\text{ if } \delta_0 \leq u \leq 1-\delta_1,\\
			1 + r_1,
				\qquad&\text{ if } u \geq 1-\delta_1.
\end{cases}
\ee
Also $r_1 \to - f'(1)>0$ as $\delta_1 \to 0$. 
If, additionally, $\chi \in [0,1)$, then we have
\be\label{e.c060503}
	\Rd(u)
		\ge \dfrac{2}{\log^2u} - \dfrac{C \log\log\sfrac1u}{\log^3\sfrac1u},
				\qquad \text{ if } u \leq \delta_0.
\ee
The constant $C$ depends only on $\chi$ and $n$.  The constants $r_0$ and $r_1$ depend on $\chi$, $n$, $\delta_0$, and $\delta_1$.
\end{lem}
Let us make two comments.  First, the term $2/\log^2u$ in (\ref{e.c060503}) is crucial for the coefficient $\sfrac52$ in the phantom front location
\be\label{mar2912}
	m_w(t)=2t-\farc{5}{2}\log t
\ee
that appears in (\ref{mar2604}) in the pulled case.   Second, the form~\eqref{e.HR} of $f$  is mainly used to prove the bound~\eqref{e.c051604}.  Indeed, the estimate~\eqref{e.c060503} follows directly from \Cref{l.eta_pulled} and the definition~\eqref{e.Q} of~$Q$.
The proof of \Cref{l.R_bound} is found in \Cref{s.Q}.

\subsection{The pushed case: the proof of \Cref{t.sdf-hr}(iii)}\label{subsec:pushed-proof}

We begin with the pushed case as it is simplest.  From~\eqref{e.sdf}, \Cref{l.R_bound}, and~\eqref{e.w>0}, we find
\be
	w_t - w_{xx} \leq w.
\ee
Hence, $e^{-t}w$ is a subsolution of the heat equation and we find, by~\eqref{e.u_0},
\be
w(t,x)
\leq e^t \int_{-\infty}^\infty w_0(x-y) \frac{e^{-\frac{y^2}{4t}}}{\sqrt{4\pi t}} dy.
\ee
As $u_0(x)=0$ for $x\ge L_0$, we also have $w_0(x)=0$ for $x\ge L_0$, and we can assume without loss of generality that
$L_0=0$. We obtain, for $x\ge 0$
\be
\bal
w(t,x)
&\leq e^t \int_{x}^\infty w_0(x-y) \frac{e^{-\frac{y^2}{4t}}}{\sqrt{4\pi t}} dy=
e^t \int_{x}^\infty \big[-u_0'(x-y)-\eta(u_0(x-y))\big] \frac{e^{-\frac{y^2}{4t}}}{\sqrt{4\pi t}} dy\\
&\le e^t \int_{x}^\infty \big[-u_0'(x-y)\big] \frac{e^{-\frac{y^2}{4t}}}{\sqrt{4\pi t}} dy=
e^t \int_{x}^\infty \partial_y\big[u_0(x-y)\big] \frac{e^{-\frac{y^2}{4t}}}{\sqrt{4\pi t}} dy\\
&= e^t \int_{x}^\infty \ u_0(x-y)\farc{y}{2t}\frac{e^{-\frac{y^2}{4t}}}{\sqrt{4\pi t}} dy\
\leq  \frac{e^t}{\sqrt{4\pi t}}  e^{- \frac{x^2}{4t}}.
\enbal
\ee
The result follows by changing variables $x \mapsto x + m(t)=x+c_*t+x_0$.~$\Box$

\subsection{The pushmi-pullyu case: the proof of \Cref{t.sdf-hr}(ii)}

We begin with the pushmi-pullyu case $\chi=1$. In that case, the front location is
\be
m(t)=2t-\farc12\log t.
\ee
We recall the following estimate to the right of $m(t)$ when $\chi=1$.
\begin{lem}\label{l.w_bound_pp}
	For any $t$ sufficiently large and any $L$, we have
	\be
		w(t,x + m(t) - L) \leq \frac{C_L}{t} (x_+ +1) e^{-x_+ - \frac{x_+^2}{Ct}}.
	\ee
\end{lem}
We omit this proof as it is essentially the same as \cite[Lemma~6.6]{AHR2}. 
In view of \Cref{l.w_bound_pp}, we need only consider the behavior of $w(t,x)$ behind the position $m(t) - L$.  
We do this via the construction of a super-solution.  Changing to the moving frame
\be
	\widetilde w(t,x) = w(t, x + m(t) - L)
		\quad\text{ and }\quad
	\tilde u(t,x) = u(t, x + m(t) - L),
\ee
and applying \Cref{l.R_bound} to (\ref{e.sdf}), we find, for any $\eps>0$,
\be
	\widetilde w_t - \Big(2 - \frac{1}{2t}\Big) \widetilde w_x
		\leq \widetilde w_{xx} + (f'(1) + \eps)\widetilde w
	\qquad\text{ for } x < 0.
\ee
Above we have potentially increased $L$ so that, by \Cref{p.weak_front}, $u > 1 -\delta_1$ with $\delta_1$ as in \Cref{l.R_bound} for $x < 0$.
 
We next remove an integrating factor.  Let $\lambda_{1,\eps}$ be the positive root of
\be\label{e.c071701}
	- 2\lambda = \lambda^2 + f'(1) + 2\eps
\ee
(cf.~\eqref{e.lambda_bis}), and let
\be
	z(t,x) = e^{- \lambda_{1,\eps}x} \widetilde w(t,x),
\ee
we obtain the differential inequality
\be\label{e.c60201}
	z_t - \Big(2(1+\lambda_{1,\eps}) - \frac{1}{2t}\Big) z_x
		\leq z_{xx}- \frac{\lambda_{1,\eps}}{2t} z - \eps z
	\qquad\text{ for } x < 0.
\ee

Before constructing a supersolution for~\eqref{e.c60201}, we note the following boundary conditions.  First, due to \Cref{l.w_bound_pp}, we have
\be
	\widetilde w(t,0)
		\leq \frac{C_{L,\eps}}{t}.
\ee
Second, due to \Cref{l.steepness} and parabolic regularity theory, we have, for any $x<0$,
\be\label{e.c051606}
	\widetilde w(t,x)
		= - \tilde u_x(t,x)
			- \eta(\tilde u(t,x))
		\leq C \sup_{(s,x) \in [t-1,t]\times [x-1, x + 1]} (1 - \tilde u(s,x))
		\leq C e^{\lambda_1 x }.
\ee
As a result, if we can produce a supersolution $\overline z(t,x)$ for~\eqref{e.c60201} defined for $t\ge T$ and $x\in [-\delta t, 0]$ that satisfies
the boundary conditions 
\be\label{e.c60203}
	\overline z(t,0)
		\geq \frac{C}{t}
	\quad\hbox{and}\quad
	\overline z(t,-\delta t)
		\geq C e^{-(\lambda_1 - \lambda_{1,\eps}) \delta t},~~\hbox{for $t\ge T$},
\ee
and the initial condition at $t=T$
\be\label{mar2816}
	 \inf_{x \in [-\delta T,0]} \overline z(t,x) \geq  C,
\ee
then we would conclude, via the comparison principle, that $\tilde z(t,x) \leq \overline z(t,x)$ for $t\ge T$ and $x\in [-\delta t, 0]$.  Let us note that $\lambda_{1,\eps} < \lambda_1$ due to~\eqref{e.c071701}.

We define the function $\overline z(t,x)$ by
\be
\overline z(t,x)= \frac{A}{t}
\qquad \text{ for } x < 0 \text{ and } t > T.
\ee
It is clearly possible to choose $A$, depending on $L$, $\delta$ and $T>0$, 
so that the conditions in~\eqref{e.c60203}-(\ref{mar2816}) are satisfied.  
It remains to check that $\overline z$ is a super-solution of~\eqref{e.c60201}.
A direct computation yields, for any $x \in (-\delta t, 0)$,
\be
\bal
\overline z_t - &\Big(2(1+\lambda_{1,\eps}) - \frac{1}{2t}\Big) \overline z_x
- \overline z_{xx} +  \Big(\frac{\lambda_{1,\eps}}{2t}+\eps\Big)\overline z
		= \overline z \left(
				- \frac{1}{t}
				+ \frac{\lambda_{1,\eps}}{2t} + \eps
				\right)
			>0,
\enbal
\ee
as long as we 
increase $T$ if necessary. 
Hence, $\overline z$ is a super-solution for~\eqref{e.c60201}. We deduce that
\be
\tilde w(t,x)\le \frac{ A}{t}e^{\lambda_{1,\eps} x},~~\hbox{ for $t>T$ and $-\delta t\le x\le 0$.}
\ee
In view of~\eqref{e.c071701}, $\lambda_{1,\eps} \nearrow \lambda_1$ as $\eps \to 0$.  Hence, the above is the desired bound for $x \in [-\delta t, 0]$. 
On the other hand, the bounds on $\tilde w$ for $x \leq -\delta t$ follow directly from~\eqref{e.c051606}.  
This completes the proof of \Cref{t.sdf-hr}(ii).~$\Box$

\subsection{The pulled case: the proof of \Cref{t.sdf-hr}(i)}

When $0\le\chi<1$ the front is located at the position 
\be\label{mar2914}
m(t)=2t-\farc{3}{2}\log t.
\ee
Exactly the same argument as in  the proof of \Cref{t.sdf-hr}(ii) to control the behavior of $w(t,x)$
for~$x <m(t)$ can be applied. Thus, we only
need to control $w(t,x+m(t))$ for $x>0$. This is done by the following. 
\begin{lem}\label{l.w_pulled}
Under the assumptions of \Cref{t.sdf-hr}(i), we have
\be
	w(t,x + m(t))
		\leq \frac{C (x^2 + 1)}{t} e^{- x - \frac{x^2}{Ct}}
		\qquad\text{ for all } x > 0.
\ee
\end{lem}
Before starting the proof, let us make the following comment. As discussed in the introduction, the convergence rate 
of $w(t,x)$ is controlled by the lag $D(t)$ of the phantom front $m_w(t)$  
behind the true front  $m(t)$, as in (\ref{jun2012})-(\ref{mar2602}). 
When $0\le\chi<1$, the phantom front $m_w(t)$ is given by~(\ref{mar2912}) 
and~$m(t)$ in (\ref{mar2914}).  On the other hand, the use of the naive linearization such as (\ref{mar2631})
\be
w_t\approx w_{xx}+w,
\ee
would produce an incorrect estimate $m_w(t)\sim 2t-(3/2)\log t$ which would lead to
$D(t)\sim O(1)$, and a bound in the spirit of (\ref{mar2635})  on the convergence rate  would be useless. 
 Thus, the lag  
comes solely from the non-zero term~$R(u)$ in~\eqref{e.c060503}. 
We have to use this estimate in an essential way to obtain any convergence rate in (\ref{e.rate}) in the pulled case, 
let alone a sharp one. 

 {\bf Proof.} First, for $L$ and $T>0$ to be determined, we let
\be\label{jun2102}
	\widetilde w(t,x)
		= w(t,x + m(t) - L)
		= w(t, x + 2t - \tfrac{3}{2}\log(t+T) - L),
\ee
and define $\tilde u$ similarly.  Then, recalling \Cref{l.R_bound}, since $\eta''(u)\le 0$, we find
\be\label{e.c60305}
	\widetilde w_t - \Big(2 - \frac{3}{2(t+T)} \Big) \widetilde w_x
		\leq \widetilde w_{xx} + (1-\Rd(\tilde u)) \widetilde w.
\ee
We remove an exponential,
\be\label{jun2016}
	z(t,x) = e^{x} \widetilde w(t,x)
\ee
to obtain
\be\label{e.c60304}
	z_t
		+ \frac{3}{2(t+T)}(z_x - z)
		\leq z_{xx} - z \Rd(\tilde u).
\ee

We now define a supersolution  to~\eqref{e.c60305} for $t\ge 1$ and $x\in\Rm$
as follows. 
For $B\geq 1$ and $T\ge 1$ to be chosen, let 
\be\label{jun2014}
	\zeta(t,x)
		= \theta(t) \Big(\frac{x+B}{B}\Big)^2
			\exp\left\{4 - 2 \sqrt{\theta(t)} - \frac{(x+B)^2}{4(t+T)}\left(1 - \frac{1}{8} \sqrt{\theta(t)}\right)\right\},
\ee 
where we have defined
\be
	\theta(t) = \frac{T}{t+T}.
\ee
Let us set
\be
	\overline w(t,x)
		= \begin{cases}
			\theta(t)
			\qquad&\text{ if } x \leq 1,\\
			\min\{\theta(t), e^{-x}\zeta(t,x)\}
				\qquad&\text{ if } x \geq 1.
		\end{cases}
\ee
The proof of Lemma~\ref{l.w_pulled} will be 
finished if we show that $w(t,x) \leq A\overline w(t,x)$, with some~$A>0$.

Before we proceed, let us  explain where (\ref{jun2014})  comes from.  
First, from~\eqref{e.c051701}, we expect $w$ to ``look like'' the solution of
\be\label{jun2018} 
\phi_t = \phi_{xx} + \phi \Big(1 - \frac{2}{\log^2\sfrac1\phi}\Big).
\ee
The traveling wave solution of this equation has the asymptotics $x^2 e^{-x}$ as $x\to+\infty$~\cite{Bou-Hen}, 
which motivates  a multiplicative factor $x^2$ in (\ref{jun2014}), 
as we have already removed an exponential factor in~(\ref{jun2016}).  
On the other hand, ``far to the right,'' we should have a Gaussian behavior, 
which motivates the $\exp\{-\sfrac{x^2}{4t}\}$ type term in (\ref{jun2014}).  
In addition,  as we have mentioned above, we expect the phantom front location $m_w(t)$ to be near the front 
location for (\ref{jun2018}), which is known to be at the position given by (\ref{mar2912}).  Thus, the lag between the true and the phantom 
fronts is $D(t)\sim \log t$. Because of that, we expect $w \sim O(\sfrac1t)$. 
This explains the  multiplicative factor $\theta(t)$ in (\ref{jun2014}).  The other terms in (\ref{jun2014}) are simply technical; in particular, the~$B$ and $T$ factors  
allow to verify the supersolution 
condition and to ``fit'' $\bar w$ above $w$ initially.

By the comparison principle applied to the differential linear inequality (\ref{e.c60304})
for $z(t,x)$, 
we will have shown that 
\be\label{mar2820}
w(t,x) \leq A\overline w(t,x),~~\hbox{for $t\ge 1$ and $x\in\Rm$,}
\ee
with some $A>0$,  if we show the following:   \\
(i) the initial comparison holds:
\be\label{mar2821}
\hbox{$w(1,x) \leq A \overline w(1,x)$ for all $x\in\R$,}
\ee
(ii) the function $\overline w(t,x)$ has the form
\be\label{mar2822}
\hbox{$\overline w(t,x) = e^{-x} \zeta(t,x)$ for $t\ge 1$ and $x\ge 10$,}
\ee
or, equivalently, we have $\theta(t)\ge e^{-x} \zeta(t,x)$ in the above region,
\\
(iii) at $x=1$ we have the opposite comparison
\be\label{mar2823}
\hbox{$e^{-1} \zeta(t,1) \geq \theta(t)$ for all $t \geq 1$,}
\ee
(iv)  the function $\theta(t)$ is a super-solution to~\eqref{e.c60305} 
for $t\ge 1$ and $x\le 10$,  and  \\
(v) the function~$\zeta(t,x)$ is a super-solution to~\eqref{e.c60304} 
for $t\ge 1$ and $x\ge 1$. 
%
%

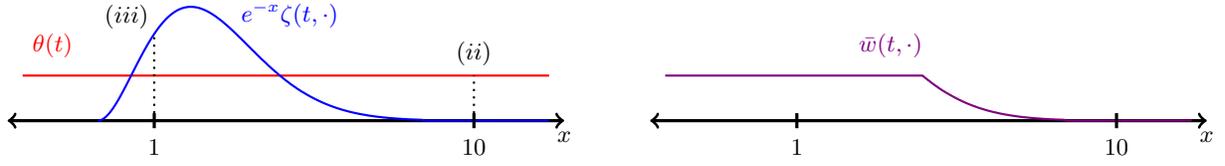
\begin{figure}  
\centering
\begin{tikzpicture}[domain=-1:6, samples=200]

  \draw[<->, very thick] (-1.2,0) -- (6.2,0) node[below] {\footnotesize$x$};

  \draw[color=red, domain = -1:6, thick]    plot (\x,.6); 
  \draw[color=blue, domain = 0:6, thick]   plot (\x,{5*\x^2 * exp(-\x - \x^2/4)})  ;
  \draw[very thick] (.75,.1)--(.75,-.1) node[below] {\footnotesize$1$};
  \draw[very thick] (5,.1)--(5,-.1) node[below] {\footnotesize$10$};
  \draw[thick, red] (-.6,.7) node[above] {\footnotesize$\theta(t)$};
  \draw[thick, blue] (2.55,1.1) node[above] {\footnotesize$e^{-x} \zeta(t,\cdot)$};
  
  \draw[thick, dotted] (.75,0) -- (.75,1.15) node[above left=-2pt] {\footnotesize$(iii)$};
  \draw[thick, dotted] (5,0) -- (5,.6) node[above] {\footnotesize$(ii)$};
\end{tikzpicture}
\qquad
\begin{tikzpicture}[domain=-1:6, samples=200]

  \draw[<->, very thick] (-1.2,0) -- (6.2,0) node[below] {\footnotesize$x$};

  \draw[color=violet, domain = -1:2.426, thick]    plot (\x,.6); 
  \draw[color=violet, domain = 2.418:6, thick]   plot (\x,{5*\x^2 * exp(-\x - \x^2/4)})  ;
  \draw[very thick] (.75,.1)--(.75,-.1) node[below] {\footnotesize$1$};
  \draw[very thick] (5,.1)--(5,-.1) node[below] {\footnotesize$10$};
  \draw[thick, violet] (2,.7) node[above] {\footnotesize$\bar w(t,\cdot)$};
    
\end{tikzpicture}
\caption{A depiction of the conditions (ii) and (iii) and their relationship to $\bar w$.}
\label{fig-shm}
\end{figure}

In particular, (\ref{mar2822})-(\ref{mar2823}) are important because they allow us to make the matching between~$\theta(t)$ and $e^{-x}\zeta(t,x)$ somewhere in 
the interval $(1,10)$  as the minimum of two super-solutions.  This is crucial because, as $\zeta(t,x)$ vanishes at $x=-B$, it can not be a super-solution for $x < 0$,
and, as we will see, $\theta(t)$ is not a super-solution for $x > 10$. This is depicted in  Figure~\ref{fig-shm}.

We now check conditions (i)-(v).  The initial comparison (\ref{mar2821}) is easy to check using well-known bounds on parabolic equations.  
In particular, $w(t,x)$ is bounded, up to a large multiplicative constant, by a Gaussian in $x$, for each $t>0$ fixed.  
Hence, after increasing $T$, independent of all parameters, and increasing $A$, depending on $L$ and $B$, the bound  (\ref{mar2821}) must hold.
Recall that $L$ appears in the change of variables (\ref{jun2102}). 

Next, we notice that (ii) is clear by observation if $B$ is sufficiently large.  
Similarly, after increasing~$T$ (depending only on $B$), (iii) is also clear by observation.  

To see that (iv) is satisfied requires us to increase $L$ (independent of all parameters) and apply \Cref{p.weak_front} with any $\delta_1$ sufficiently small to find that
\be
	\tilde u(t,x) \geq 1- \delta_1
		\qquad\text{ for all } t \geq 1,~x \leq 10.
\ee
Then, from~\Cref{l.R_bound}, we have
\be
	1- \Rd(\tilde u) \leq -r_1
		\qquad\text{ for all } t \geq 1, ~x \leq 10.
\ee
Thus, up to increasing $T$, depending only on $\delta_1>0$, we have
\be
	\theta_t
		- \Big(2 - \frac{3}{2(t+T)}\Big) \theta_x
		- \theta_{xx}
		- (1 - \Rd(\tilde u)) \theta
	\ge -\frac{T}{(t+T)^2}
			+ r_1\frac{T}{t+T}
		> 0.
\ee
Therefore, (iv) holds.

We now check (v), which is a computationally tedious condition to verify, even though the computations are completely elementary.   
First, we compute:  
\be
\begin{split}
		&\frac{\overline z_t + \frac{3}{2(t+T)}(\overline z_x - \overline z)
			- \overline z_{xx} + \overline z \Rd(\tilde u)}{\overline z}
			 = \frac{\dot \theta}{\theta}
				- \frac{\dot\theta}{\sqrt \theta}
				+ \frac{(x+B)^2}{4(t+T)^2} \big(1 - \frac{1}{8}\sqrt\theta\big)
				+ \frac{1}{8} \frac{\dot \theta}{2\sqrt \theta} \frac{(x+B)^2}{4(t+T)}
			\\&\qquad\quad
				+ \frac{3}{2(t+T)} \Big(\frac{2}{x+B} - \frac{x+B}{2(t+T)}\big(1 - \frac{1}{8} \sqrt \theta\big) - 1\Big)
			\\&\qquad\quad
				- \left(\frac{2}{(x+B)^2} - \frac{5}{2(t+T)} \big(1 - \frac{1}{8} \sqrt \theta\big) +  \frac{(x+B)^2}{4(t+T)^2}\big(1 - \frac{1}{8}\sqrt \theta\big)^2\right) + \Rd(\tilde u).
	\end{split}
\ee
Noticing that $\dot\theta/\theta = -1/(t+T)$ and $\dot\theta/\sqrt\theta = -\sqrt\theta/(t+T)$, cancelling the obvious terms, and then grouping terms by the growth in $x$ yields
\be
	\begin{split}
		&\frac{\overline z_t - \frac{3}{2(t+T)}(\overline z_x - \overline z)
			- \overline z_{xx} + \overline z \Rd(\tilde u)}{\overline z}
			\\&\quad
			 =
				\Big(\frac{\sqrt\theta}{t+T}- \frac{3}{2(t+T)} \frac{1}{8} \sqrt \theta\Big)
				+ \frac{(x+B)^2}{4(t+T)^2} \Big(\big(1 - \frac{1}{8}\sqrt\theta\big)
				- \frac{1}{2(t+T)} \Big) \frac{1}{8} \sqrt \theta
			\\&\qquad\quad
				+ \frac{3}{2(t+T)} \Big(\frac{2}{x+B} - \frac{x+B}{2(t+T)}\big(1 - \frac{1}{8} \sqrt \theta\big)\Big)
				- \frac{2}{(x+B)^2}
				 + \Rd(\tilde u).
	\end{split}
\ee
Since $\theta \leq 1$, we have, up to increasing $T$ (independent of all parameters),
\be
	\begin{split}
		&\frac{\overline z_t - \frac{3}{2(t+T)}(\overline z_x - \overline z)
			- \overline z_{xx} + \overline z \Rd(\tilde u)}{\overline z}
			\\&\qquad
			 \geq
				\frac{\sqrt\theta}{2(t+T)}
				+ \frac{(x+B)^2}{4(t+T)^2}  \frac{\sqrt \theta}{16}
				+ \frac{3}{2(t+T)} \Big(\frac{2}{x+B} - \frac{x+B}{2(t+T)}\Big)
				- \frac{2}{(x+B)^2}
				 + \Rd(\tilde u).
	\end{split}
\ee
Using Young's inequality and then increasing $T$ (independent of all parameters), we arrive at
\be\label{e.c60401}
	\begin{split}
		&\frac{\overline z_t - \frac{3}{2(t+T)}(\overline z_x - \overline z)
			- \overline z_{xx} + \overline z \Rd(\tilde u)}{\overline z}
			\\&\qquad
			 \geq
				\frac{\sqrt\theta}{2(t+T)}
				+ \frac{(x+B)^2}{8(t+T)^2}  \frac{\sqrt \theta}{16}
				+ \frac{3}{2(t+T)}\frac{2}{x+B}
				- \frac{72 \sqrt \theta}{T (t+T)}
				- \frac{2}{(x+B)^2}
				 + \Rd(\tilde u)
			\\&\qquad
			 \geq
				\frac{\sqrt\theta}{4(t+T)}
				+ \frac{(x+B)^2}{8(t+T)^2}  \frac{\sqrt \theta}{16}
				+ \frac{3}{2(t+T)}\frac{2}{x+B}
				- \frac{2}{(x+B)^2}
				 + \Rd(\tilde u).
	\end{split}
\ee

At this point, we can see why the right hand side of~\eqref{e.c60401} should be positive.  
Recall that, according to \Cref{l.R_bound} (equation~\eqref{mar2812}), the term $\Rd(\tilde u) \ge r_0>0$ 
when $\tilde u$ is not too small.
Hence, it should dominate the next to last term in the right side of (\ref{e.c60401}) in that region if $B$ is large. 
On the other hand, for $\tilde u$ small, the term $\Rd(\tilde u)$ looks like $2/\log^2(\tilde u)$, according to (\ref{e.c060503}).  
Moreover, as~$\tilde u(t,x) \approx U_*(x)$ and $U_*(x)$ has the asymptotics~\eqref{e.lambda}, we have $\log^2(\tilde u) \approx x^2$.
Thus, once again,~$R(\tilde u)$ dominates the next to last term in the right side of (\ref{e.c60401}).

We make the discussion above more precise.  Let us  fix $\delta_1>0$  as in \Cref{l.R_bound}.  
We claim that, up to increasing $L$ (depending on $\delta_1$), we have 
\be\label{e.c60402}
	\tilde u(t,x)
		\geq \begin{cases}
				1-\delta_1
					\qquad&\text{ if } x \leq L/2,\\
				\frac{x + 1}{C_L} e^{-x - \frac{C_Lx^2}{t}}
					\qquad&\text{ if } x \geq L/2,
			\end{cases}
\ee
for all $t\ge 1$, with a constant $C_L$ that depends on $L$.  The first alternative above is due to \Cref{p.weak_front}.  
The second alternative follows from \cite[Proposition~3.1]{HNRR} and its proof, as well as an application of the comparison principle.

We first consider the ``large'' $\tilde u$ regime (and, thus, $x$ ``not too far on the right'').  If $\tilde u \geq \delta_0$, then~$\Rd(\tilde u) \geq r_0$ 
due to (\ref{mar2812}) and we find
\be
	-\frac{2}{(x+B)^2} + \Rd(\tilde u)
		\geq -\frac{2}{(x+B)^2} + r_0
		> 0,
\ee
up to increasing $B$ further if necessary so that $2/B^2 < r_0$.  
In particular, then we have, from~\eqref{e.c60401},
\be
	\frac{\overline z_t - \frac{3}{2(t+T)}(\overline z_x - \overline z)
			- \overline z_{xx} + \overline z \Rd(\tilde u)}{\overline z}
			> 0,~~\hbox{ if $\tilde u(t,x)\ge\delta_0$,} 
\ee
as desired.  

Next we consider the ``small'' $\tilde u$ regime (and, thus, ``large'' $x$ regime). Note that, by~\eqref{e.c60402}, if~$\tilde u \leq \delta_0$, then
\be\label{e.c032901}
	x \geq \min\Big( \frac{1}{2} \log\frac{1}{C\delta_0}, \sqrt{\frac{t}{2C} \log\frac{1}{C\delta_0}}\Big)
		\geq \sqrt{\frac{1}{2C} \log\frac{1}{C\delta_0}}.
\ee
In particular, this case is restricted to $x$ that is very large,  after possibly decreasing $\delta_0$.

We begin by estimating $\Rd(\tilde u)$ using~\eqref{e.c060503}.  
For the  quadratic term, we apply~\eqref{e.c60402} to find
\be
	\frac{2}{(\log(\tilde u))^2}
		\geq \frac{2}{x^2} \frac{1}{(1 + \frac{Cx}{t} - \frac{\log x}{x} + \frac{\log C}{x})^2}.
\ee
Then, using that $(1+z)^{-2} \geq 1- 2z$ for all $z\geq -1$, we obtain
\be
\frac{2}{(\log(\tilde u))^2}
		\geq \frac{2}{x^2}\left(1 - 2\left(\frac{Cx}{t} - \frac{\log x}{x} + \frac{\log C }{x}\right)\right)
		= \frac{2}{x^2} - \frac{4C}{x t} + \frac{4 \log x}{x^3} - \frac{4 \log C}{x^3}.
\ee
A similar argument, using the inequality 
\be
(1-z)^{-3} \leq 1 + Cz,~~ \hbox{for $0\le z \leq 1/2$}, 
\ee
yields a bound for the second term in $R(\tilde u)$:
\be
	\begin{split}
		- \frac{C}{|\log(\tilde u)|^3}
			&\geq -\frac{C}{x^3} \frac{1}{\left(1 + \frac{Cx}{t} - \frac{\log x}{x} + \frac{\log C}{x} \right)^3}
			\geq -\frac{C}{x^3} \frac{1}{\left(1 - \frac{\log x}{x}\right)^3}
			\\&
			\geq - \frac{C}{x^3}\left(1 +C \frac{\log x}{x} \right)
			= - \frac{C}{x^3} - \frac{C \log x}{x^4} .
	\end{split}
\ee
Using these in~\eqref{e.c60401}, we find
\be
	\begin{split}
		&\frac{\overline z_t - \frac{3}{2(t+T)}(\overline z_x - \overline z)
			- \overline z_{xx} - \overline z \Rd(\tilde u)}{\overline z}
			 \geq
				\frac{\sqrt\theta}{4(t+T)}
				+ \frac{(x+B)^2}{8(t+T)^2}  \frac{\sqrt \theta}{16}
				+ \frac{3}{(t+T)(x+B)}
				\\&
				\phantom{MMMMM}
				- \frac{2}{(x+B)^2}
				+ \frac{2}{x^2} - \frac{4C}{x t} + \frac{4\log x}{x^3}
				- \frac{4 \log C}{x^3}
				- \frac{C}{x^3} - \frac{C \log x}{x^4}. 
	\end{split}
\ee
After decreasing $\delta_0$ (which, by~\eqref{e.c032901}, increases the lower bound for $x$), we find
\[
	\begin{split}
		&\frac{\overline z_t - \frac{3}{2(t+T)}(\overline z_x - \overline z)
			- \overline z_{xx} - \overline z \Rd(\tilde u)}{\overline z}
			 \geq
				\frac{\sqrt\theta}{4(t+T)}
				+ \frac{(x+B)^2}{8(t+T)^2}  \frac{\sqrt \theta}{16}
				+ \frac{3}{(t+T)(x+B)}
				- \frac{4C}{x t} + \frac{2\log x}{x^3}.
	\end{split}
\]
There is only one negative term above.  Applying Young's inequality with $p=3/2$ and $q = 3$  yields
\be
	- \frac{4C}{x t}
		\geq -\frac{\sqrt\theta}{4(t+T)}
			- C\bigg(\Big(\frac{\sqrt\theta}{(t+T)}\Big)^{-\frac{2}{3}} \frac{1}{xt}\bigg)^3
		= -\frac{\sqrt\theta}{4(t+T)}
			- C\frac{(t+T)^3}{T} \frac{1}{x^3t^3}
		\geq -\frac{\sqrt\theta}{4(t+T)}
			- C\frac{T^2}{x^3}.
\ee
Hence, we have
\be
	\begin{split}
		&\frac{\overline z_t - \frac{3}{2(t+T)}(\overline z_x - \overline z)
			- \overline z_{xx} - \overline z \Rd(\tilde u)}{\overline z}
			 \geq
				\frac{(x+B)^2}{8(t+T)^2}  \frac{\sqrt \theta}{16}
				+ \frac{3}{2(t+T)}\frac{2}{x+B}
				- C\frac{T^2}{x^3} + \frac{\log x}{2x^3}.
	\end{split}
\ee
which is positive after further decreasing $\delta_0$ (which, by~\eqref{e.c032901}, increases $x$).   
This concludes the proof of (v) and, thus, the proof of the lemma.~$\Box$

\section{Proofs of the bounds on $\eta$ and $Q$}\label{s.Q}

\subsection{Concavity of $\eta$: \Cref{prop-mar2602}}\label{s.Q_concavity}

We make two observations.  First, arguing as in \Cref{l.eta_pulled}, it is easy to check that, for any $f$, its traveling wave profile function $\eta$ satisfies
\be\label{mar2325}
	\eta^2(u)\eta''(u)\to 0,~~\hbox{as $u\to 1^-$}.
\ee 	
Second,   \Cref{prop-mar2602} follows from the following more general result.
\begin{lem}\label{l.general_eta}
	Assume that \eqref{e.rde}-\eqref{e.normalization} hold.  Suppose that either:
	\begin{enumerate}[(i)]
		\item (pulled case) the asymptotics~\eqref{e.lambda} hold  and $f'' \leq 0$ on $(0,1)$;
		\item (pulled case) the asymptotics~\eqref{e.lambda} hold  
		and there is $u_0 \in [0,1]$ such that $f'' \geq 0$ on $(0,u_0)$ and $f'' \leq 0$ on $(u_0,1)$;
		\item (pushed and pushmi-pullyu cases) there is $\chi \geq 1$ and $A$ satisfying $A(0) = A'(0) = 0$ and $A(1) = 1$ such that
		\be\label{e.c062505}
			f(u) = (u-A(u))(1 + \chi A'(u))
			\quad\text{ and }\quad
			A'', A''' \geq 0.
		\ee
		If $\chi = 1$, the condition $A''' \geq 0$ is not necessary.
	\end{enumerate}
	Then $\eta'' \leq 0$ and $Q \leq 1$.
\end{lem}

\noindent{\bf Proof in cases (i) and (ii).} First, note that, case (i) is really the 
subcase of (ii) where $u_0 = 0$.  Hence, we only consider
case (ii).  Let us also recall that $f'(0)=1$, according to assumption (\ref{e.normalization}). If 
the asymptotics (\ref{e.lambda}) holds and $c\ge 2$ is the speed of the wave,  
then, by linearization as $x\to+\ifnty$, it is easy to see that $\lambda_0$  must be a double root of the equation 
\be
c\lambda=\lambda^2+1.
\ee
It follows that $c=c_*=2$ and $\lambda_0=1$. 
 
Observe that it is thus enough to show that $\eta'' \leq 0$.  Indeed,
\be
	Q
		= \eta'(2 - \eta') + \eta'' \eta
		\leq 1 + \eta'' \eta
		\leq 1.
\ee
By \Cref{l.eta_pulled}(iii), there exists $u_1>0$ so that
\be\label{mar2702}
\eta''(u)\le 0,~~\hbox{ for all $0<u<u_1$.}
\ee
Thus, the following is well-defined and positive:  
\be\label{e.bar_u}
\bar u= \sup\{ \tilde u\in(0,1):~\eta''(u) \leq 0 \text{ on } (0,\tilde u]\}.
\ee
Our goal is to prove that $\bar u=1$. 

Suppose, for the sake of a contradiction, that $\bar u<1$. 
Writing~\eqref{e.f_eta} as 
\be\label{e.eta_ode}
2-\eta'=\farc{f(u)}{\eta(u)},
\ee
we find $\eta^2 \eta'' = \eta' f - f' \eta$ and, hence,
\be\label{e.c031902}
(\eta^2 \eta'')' =(\eta' f-f'\eta)'= \eta'' f - \eta f''.
\ee
It follows that, at $\bar u$, we have
\be\label{e.c031901}
0 \leq (\eta^2 \eta'')' (\bar u)= \eta''(\bar u) f(\bar u) - \eta(\bar u) f''(\bar u)
= - \eta(\bar u)f''(\bar u).
\ee
The first inequality follows from the fact that $\eta^2 \eta''$ crosses zero at $\bar u$ due to~\eqref{e.bar_u}. 
As $\eta(\bar u) > 0$ (recall~\eqref{mar2618}), it follows that $f''(\bar u)\le 0$, which in turn implies that 
\be\label{mar2627}
\bar u \geq u_0.
\ee  
We deduce that
\be\label{mar2706}
f''(u)\le 0~~\hbox{ for all $u\ge\bar u$.}
\ee

We now claim that $\eta'' > 0$ on $(\bar u, 1)$.  The definition \eqref{e.bar_u} of $\bar u$ implies that if $\bar u<1$ then for every $\eps>0$
sufficiently small, there is~$u_\eps \in (\bar u, \bar u + \eps)$ such that~$\eta''(u_\eps)>0$. 
Suppose that there is~$\bar v_\eps\in(u_\eps,1)$ such that~$\eta''(u)>0$ for~$u\in(u_\eps,\bar v_\eps)$ and~$\eta''(\bar v_\eps)=0$. Then, integrating (\ref{e.c031902})
gives
\be\label{mar2624}
	0 > - \eta^2 (u_\eps)\eta''(u_\eps)
		= \int_{u_\eps}^{v_\eps} \left( \eta'' f - \eta f'' \right) du
		> 0,
\ee
which is a contradiction. 
The second inequality in~\eqref{mar2624} follows from the fact that, on the domain on integration, $\eta, f, \eta'' > 0$ and $f'' \leq 0$.   
 We conclude that~$\eta''(u) > 0$ for $u \in (u_\eps, 1)$.  
By the arbitrariness of $\eps>0$, it follows that $\eta'' > 0$ on~$(\bar u,1)$, as claimed.

Finally, we conclude by obtaining a contradiction at $u=1$.  Going back to (\ref{e.c031902}) and recalling~(\ref{mar2627}), we deduce that 
\be\label{mar2628}
(\eta^2 \eta'')' =\eta'' f - \eta f''>0,~~\hbox{for $\bar u< u< 1$.}
\ee
Recall that $\eta''(\bar u) = 0$, by construction.  As a consequence, we obtain, for any $u > \bar u$,
\be\label{mar2628bis}
	\eta^2(u)\eta''(u)
		= \eta^2(\bar u)\eta''(\bar u) + \int_{\bar u}^u \left( \eta'' f - \eta f''\right) du
		= \int_{\bar u}^u \left( \eta'' f - \eta f''\right) du
		> 0. 
\ee
Taking the limit $u \nearrow 1$ and using~\eqref{mar2325}, we obtain
\be
	0
		= \lim_{u\nearrow 1} \eta^2(u) \eta''(u)
		= \int_{\bar u}^1 \left( \eta'' f - \eta f''\right) du
		>0.
\ee
Here, the last inequality follows from~\eqref{mar2628}. 
This  contradiction shows that it is impossible that~$\bar u<1$.  
It follows that $\bar u=1$ and  $\eta''(u)<0$ for all $u\in(0,1)$.  This concludes the proof.~$\Box$

\bigskip

\noindent{\bf Proof in case (iii).} Here, we have the explicit form of $\eta$ due to \cite[Proposition~A.2]{AHR2}:
\be\label{e.c060501}
	\eta(u) = \sqrt \chi (u - A(u))
	\quad\text{ and }\quad
	c_* = \sqrt \chi + \frac{1}{\sqrt \chi}.
\ee
It is immediate that $\eta'' \leq 0$; hence, we need only show that $Q \leq 1$.  A direct computation yields
\be
	Q
		= 1 + (\chi - 1)A' - \chi |A'|^2 - \chi (u-A) A''
		\leq 1 + \chi( A' - |A'|^2 - (u-A)A'').
\ee
The second inequality follows from the convexity of $A$ and the fact that $A'(0) = 0$, which imply that~$A'\geq 0$.  It is, hence, enough to show that 
\be\label{jun2106}
A' - |A'|^2 - (u-A)A'' \leq 0.
\ee
Note that, at $u=0$, the expression above vanishes  On the other hand,
\be\label{jun2202}
	(A' - |A'|^2 - (u-A)A'')'
		= - A' A'' - (u - A) A''' \leq 0,
\ee
since $u-A, A', A'', A''' \geq 0$.  We conclude that $Q\leq 1$. This completes the proof.

Finally, we consider the last statement for $\chi = 1$.  We have already observed that $\eta'' \leq 0$.   We conclude by noting that, from~\eqref{e.c060501}, $c_* = 2$ and then arguing as in the second paragraph of the proof for cases (i) and (ii).~$\Box$


\subsection{Refined bounds on $Q$: proof of \Cref{l.R_bound}}

First, we note that the bounds in~\eqref{e.c051604} follow from \Cref{l.general_eta}.  Second, the bounds~\eqref{e.c060503} follow directly from \Cref{l.eta_pulled}.

We now address the bounds in~\eqref{mar2812} for the remainder of the proof.  We first investigate the first alternative in~\eqref{mar2812}.  In the case $\chi \geq 1$, the proof of \Cref{l.general_eta} clearly shows 
that if $A'', A''' < 0$, then $Q$ is bounded away from $1$ on compact subsets of $(0,1]$.  This is exactly the first alternative in~\eqref{mar2812} for the case $\chi \geq 1$.

When $0 \leq \chi < 1$, the first inequality in~\eqref{mar2812} is deduced using only the concavity of $\eta$~\eqref{e.c051604} and the asymptotics \Cref{l.eta_pulled}.  Indeed, these imply that $\eta'(u) \leq \eta'(\delta_0) < 1$ for all $u \in (\delta_0,1)$.  Hence,
\be\label{e.c062702}
	\begin{split}
		R(u) &= 1 - Q(u)
			= 1 - \eta'(u)(2 - \eta'(u)) - \eta(u) \eta''(u)
			\geq 1 - \eta'(u)(2 - \eta'(u))
			\\&
			> 1 - \eta'((\delta_0)(2 - \eta'(\delta_0))
			> 0.
	\end{split}
\ee
This yields the first alternative in~\eqref{mar2812} in the pulled case.

We now investigate the second alternative in~\eqref{mar2812}.  Notice that
\be\label{e.c060701}
		Q(1) = \eta'(1) (c_* - \eta'(1)) + \eta(1)\eta''(1)
			= -\lambda_1 (c_* + \lambda_1) 
			= f'(1) < 0.
\ee
The second equality above follows from~\eqref{e.eta} and~\eqref{mar2614}, while the third is due to~\eqref{e.lambda_bis}.  The inequality uses the particular form of $f$.  
This concludes the proof.~$\Box$

\section{The general case}\label{s.general}

In this section, we discuss the convergence rate when 
$f$ satisfies~\eqref{e.f} and the normalization (\ref{e.normalization}) but does not necessarily have the Hadeler-Rothe form (\ref{e.HR}).

Let us begin by recalling the proof of the convergence rates in the Hadeler-Rothe case (\Cref{t.rate-hr}).  The main lemma is the estimate on $w$ (\Cref{t.sdf-hr}).  The argument to deduce \Cref{t.rate-hr} from \Cref{t.sdf-hr} relies only on the behavior of $\eta$ near $u=0$, which is established in full generality in \Cref{l.eta_pulled,lem.eta_pushed}.

In the proof of \Cref{t.sdf-hr}, there are exactly two places where we use the assumption (\ref{e.HR}) on the form of $f$  
rather than just the assumptions~(\ref{e.f})-(\ref{e.normalization}): 
the $O(1)$ asymptotics for the front location of $u(t,x)$ (\Cref{p.weak_front}) and the bounds on $Q$ (\Cref{l.R_bound}).  The final conclusion of \Cref{l.R_bound}, that is, the expansion~\eqref{e.c060503}, holds for any pulled front as it merely reflects the linear factor in~\eqref{e.c062202}.  Hence, the supersolutions for $w$ constructed 
in each case in the proof of \Cref{t.sdf-hr} hold in generality if we take the front asymptotics of $u$ and behavior of $\eta$ as assumptions.  Hence, the {\em exact arguments} above yield the following:
\begin{thm}\label{t.rate-general}
	Suppose that $u$ solves~\eqref{e.rde} with $f$ satisfying~\eqref{e.f}-\eqref{e.normalization} and initial data $u_0$ satisfying~\eqref{e.u_0}. 
	 Assume further that the traveling wave profile function $\eta$ and the associated quantity $Q$, respectively defined in~\eqref{e.eta} and~\eqref{e.Q}, satisfy~\eqref{e.c051604}-\eqref{mar2812}.  (Here, we are only assuming the positivity of $r_1$, not necessarily the limiting behavior as $\delta_1 \to 0$ stated below~\eqref{mar2812}.)
	Finally, suppose the front asymptotics of $u$ are given by
	\be\label{e.c062504}
		m(t) = \begin{cases}
				2t - \frac{3}{2} \log t + O(1)
					&\qquad\text{ if  $U_*$ is pulled},\\
				2t - \frac{1}{2} \log t + O(1)
					&\qquad\text{ if  $U_*$ is pushmi-pullyu type},\\
				c_*t + O(1)
					&\qquad\text{ if  $U_*$ is pushed},\\
			\end{cases}
	\ee
	in the sense of~\eqref{mar2802}, with the definition of pushed, pulled, and pushmi-pullyu given in~\eqref{e.c062202}-\eqref{e.c062203}.  Then there is $\sigma: [0,\infty) \to \R$ such that, whenever $c_* = 2$,
	\be\label{e.c062501}
		\|u(t, \cdot + \sigma(t)) - U_*(\cdot)\|_{L^\infty} \leq \frac{C}{t},
	\ee
	and, for any $\Lambda>0$, whenever $c_* > 2$,
	\be\label{e.c062502}
		\|u(t, \cdot + \sigma(t)) - U_*(\cdot)\|_{L^\infty([-\Lambda, \infty))} \leq \frac{C_\Lambda}{\sqrt t} e^{- \frac{c_*^2 - 4}{4} t}.
	\ee
	\end{thm}

\subsection{The assumptions in \Cref{t.rate-general}}

In this section, we discuss the three main assumptions in \Cref{t.rate-general}:~\eqref{e.c062504},~\eqref{e.c051604}, and~\eqref{mar2812}.  
Briefly, the front asymptotics~\eqref{e.c062504} of $u$ is nearly known in complete generality so it is a quite weak assumption and the refined bounds~\eqref{mar2812} on $Q$ may be side-stepped by alternate arguments at the expense of a slightly less precise convergence rate.  Thus, the main assumption to be checked in practice is~\eqref{e.c051604}, that is, that $\eta'' \leq 0$ and $Q\leq 1$.  We formulate a version of \Cref{t.rate-general} that assumes only~\eqref{e.c051604}
in \Cref{t.rate-general2} below.  We also discuss here the feasibility of~\eqref{e.c051604}.

\medskip
\noindent {\bf The assumption~\eqref{e.c062504}.}
In fact,~\eqref{e.c062504} is nearly established in full generality.  The pulled and pushed asymptotics in~\eqref{e.c062504} are completely proved: see~\cite[Lemma 5.2]{Uchiyama} for the pushed case and~\cite{Giletti} for the pulled case.  The statement in~\cite{Giletti} additionally requires $f'(1) <0$, although this can likely be removed via a comparison argument with solutions to~\eqref{e.rde} with appropriately chosen $\underline f$ and $\overline f$ in place of $f$.  
We do not pursue this further here.

The pushmi-pullyu case is more delicate.  If $f$ has the particular form~\eqref{e.c062505}, this is established in~\cite{AHR2}, but it is otherwise still open.  The most general result is~\cite{Giletti}, in which the 
asymptotics 
\be\label{e.c062601}
	m(t) = 2t - \frac{1}{2} \log t + o(\log t)
\ee
is established with no assumptions on $f$ beyond~\eqref{e.f}-\eqref{e.normalization}, that $U_*$ is pushmi-pullyu type, and that~$f'(1) < 0$.  Using~\eqref{e.c062601} in our arguments, we get
\be\label{e.c062607}
	D(t) = m(t) - m_w(t)
		= \log t + o(\log t).
\ee
It is not hard to track the effect of the $o(\log t)$ term in our computations to see that the informally derived convergence rate~\eqref{mar2602} holds: for every $\eps>0$,
\be\label{e.c062703}
	\|u(t, \cdot + \sigma(t)) - U_*(\cdot)\|_{L^\infty} \leq \frac{C_\eps}{t^{1-\eps}}.
\ee
Another argument leading to~\eqref{e.c062703} is sketched in greater detail below, see~\eqref{e.c062606} and its discussion.

\medskip
\noindent {\bf The assumption~\eqref{e.c051604}.}
The main purpose of this assumption is to guarantee that~\eqref{jun2304} holds; that is,
\be\label{e.c062603}
	w_t \leq w_{xx} + w.
\ee
We were unable to obtain more general assumptions on $f$ guaranteeing~\eqref{e.c051604} than those stated in Lemma~\ref{l.general_eta}.  
Further, it is not difficult to construct nonlinearities $f$ for which~$Q \not \leq 1$ (see \Cref{s.Q_not_<_1}), although these examples appear fairly pathological. 
Numerical experiments indicate that~\eqref{e.c051604} is ``often'' true. For example, consider generalizations of the Hadeler-Rothe nonlinearities of the form 
\be
f(u)=(u-A(u))(1+\chi A'(u))
	\qquad\text{ with } \chi \geq 0, A(0)= A'(0) = 0, A(1) = 1,
\ee
analyzed in~\cite{AHR2}.  Then, numerics indicates that~\eqref{e.c051604} holds as long as~$A(u)$ is increasing, convex and~$A'''(u)\geq 0$. 
For $\chi\ge 1$ this follows immediately from \Cref{l.general_eta}. 

We expect that, in many applications, either the assumptions of Lemma~\ref{l.general_eta} would hold, or 
the inequalities in~\eqref{e.c051604} are checkable or can be sidestepped using {\em ad hoc} adjustments to our approach here.   
It is easy to derive several differential equations relating $\eta$ and $Q$ to $f$ that are useful for understanding $\eta$ and $Q$, although we do not discuss this further here.

\medskip
\noindent {\bf The assumption~\eqref{mar2812}.}
This assumption is not used in the argument of the pushed setting (see~\Cref{subsec:pushed-proof}).  Hence, we need only address the pulled and pushmi-pullyu cases.

In the pushmi-pullyu case, we outline an argument below that yields nearly the same conclusion albeit without either inequality in assumption~\eqref{mar2812}.  Hence, we focus our discussion mainly on the pulled case, where~\eqref{mar2812} plays a greater role.

In the pulled setting, the first inequality in~\eqref{mar2812} holds automatically due to the concavity of~$\eta$~\eqref{e.c051604} (see~\eqref{e.c062702} and the arguments surrounding it). 
The second inequality in~\eqref{mar2812} is equivalent to $Q(1)<0$, which holds if and only if $f'(1) < 0$ (see~\eqref{e.c060701}).  

We  outline a slightly less precise argument that proceeds without the assumption~\eqref{mar2812}.
If~\eqref{e.c051604} holds, then (\ref{e.sdf}) yields~\eqref{e.c062603}.
Consider first the pushmi-pullyu case.  
The arguments in~\cite{HNRR} readily yield that any bounded subsolution of~\eqref{e.c062603} satisfies
\be
	w(t,x + m_w(t))
		= w(t,x + 2t - \sfrac32 \log t) \leq C x e^{-x - \frac{x^2}{Ct}}
		\qquad\text{ for all } x>1,
\ee
which implies the nearly sharp bound
\be\label{e.c062606}
	w(t, x + D(t) + m_w(t))
		= w(t,x + 2t - \sfrac12 \log t)
		\leq \frac{C (x + \log t)}{t} e^{-x - \frac{x^2}{Ct}}.
\ee
It appears likely that the $\log t$ error term may be avoided by using the traveling wave ``trace-back'' arguments of~\cite{HNRR} (see the proof of Proposition~3.1 therein); however, we do not pursue that here.

Note that the above estimate is assuming that $m(t) = 2t - \sfrac12 \log t + O(1)$.  If we only have available~\eqref{e.c062601} 
the above changes only by a multiplicative factor of $t^\eps$, for any $\eps>0$, due to a spatial shift of $\eps \log t$. 

The bound~\eqref{e.c062606} matches the estimate in \Cref{t.sdf-hr}(ii) up to an extra $\log t$ multiplicative factor.   
Thus, the argument deducing \Cref{t.rate-hr} from \Cref{t.sdf-hr} proceeds in the exact same manner and gives an error bound between $u$ and $U_*$ 
of the form $O(t^{-1} \log t)$ at and beyond the front.

The pulled case follows similarly, using the work of~\cite{Bou-Hen} in place of~\cite{HNRR}, and leads to an error bound of the form $O(t^{-1} \log^2 t)$  at and beyond the front. 

From the above, we deduce the following more general, slightly less precise result.  Given that the proof follows that of \Cref{t.rate-hr} exactly up to the modifications outlined above, we omit it.
\begin{thm}\label{t.rate-general2}
	Suppose that $u$ solves~\eqref{e.rde} with $f$ satisfying~\eqref{e.f}-\eqref{e.normalization} and initial data $u_0$ satisfying~\eqref{e.u_0}. 
	 Assume further that the traveling wave profile function $\eta$ and the associated quantity $Q$, defined in~\eqref{e.eta} and~\eqref{e.Q} respectively, satisfy~\eqref{e.c051604}.    
	Then there is $\sigma: [0,\infty) \to \R$ such that, for any $\Lambda>0$, $\eps>0$, and $t\geq e$,
	\be
		\|u(t, \cdot + \sigma(t)) - U_*(\cdot)\|_{L^\infty([-\Lambda, \infty))} \leq 
			\begin{cases}
				\displaystyle\frac{C_\Lambda \log (t)}{t}
					&\quad\text{ whenever $U_*$ is pulled,}\smallskip
				\\
				\displaystyle\frac{C_{\Lambda, \eps}}{t^{1-\eps}}
					&\quad\text{ whenever $U_*$ is pushmi-pullyu type,}\smallskip
		\\
				\displaystyle \frac{C}{\sqrt t} e^{- \frac{c_*^2 - 4}{4} t}
					&\quad\text{ whenever $U_*$ is pushed}.
		\end{cases}
	\ee
	If $U_*$ is pushmi-pullyu type and we additionally assume that $m(t) = 2t - \frac{1}{2}\log t + O(1)$, then the improved estimate with rate $C_\Lambda t^{-1} \log^2 t$ holds.
	\end{thm}

\subsection{An example where $Q \not < 1$}\label{s.Q_not_<_1}

\Cref{l.R_bound} is a crucial aspect of the proof of \Cref{t.rate-hr}.  As discussed above, its main component is that $Q\leq 1$, which is true for the Hadeler-Rothe nonlinearities~\eqref{e.HR}, as well as many other nonlinearities.  We show here that it is not  true for some $f$ satisfying (\ref{e.f})-(\ref{e.normalization}).
\begin{prop}
	There is a nonlinearity $f\in C^2([0,1])$ satisfying~\eqref{e.f}-\eqref{e.normalization} such that
		\be
			\sup_{u\in(0,1)} Q(u) > 1.
		\ee
\end{prop}
\noindent{\bf Proof.} 
	Let $\eta_{\rm FKPP}$ be the traveling wave profile function associated to the classical Fisher-KPP nonlinearity $f_{\rm FKPP}(u) = u(1-u).$ 
	Note that $\eta_{\rm FKPP}$  is concave by Lemma~\ref{l.general_eta}.  
	Define
	\be
		M = \sup_{u\in[0,1]} \eta_{\rm FKPP}(u),
	\ee
	and let $0<u_1 < u_2 < u_3 < u_4<1$ be such that
	\be
		\eta_{\rm FKPP}(u) \geq \frac{2M}{3}
			~~\text{ for all } u \in [u_2, u_3]
		\quad\text{ and }\quad
		\eta_{\rm FKPP}(u) \leq \frac{M}{3}
			~~\text{ for all } u \in [0,u_1] \cup [u_4,1].
	\ee 
	This definition is illustrated in \Cref{f.eta}.

	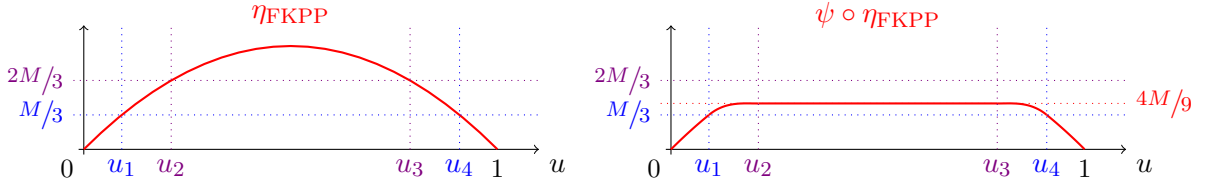
\begin{figure}
	\centering
	\begin{tikzpicture}[scale=5.5]
  \draw[->] (-0.0125,0) -- (1.1,0) node[below right] {$u$};
  \draw[->] (0,-0.0125) -- (0,0.3) node[left] {};

  \draw[red,thick, domain=0:1] plot (\x, {\x*(1-\x)});

  \draw[blue, dotted] (-0.025, 1/12) -- (1.1, 1/12);
  \draw[violet, dotted] (-0.025, 1/6) -- (1.1, 1/6);

  \draw[blue, dotted] (.0917517, 0) -- (.0917517, 0.3);
	  \node[blue, below] at (.0917517,0) { $u_1$};
  \draw[violet, dotted] (.2113249, 0) -- (.2113249, 0.3);
	  \node[violet, below] at (.2113249,0) { $u_2$};
  \draw[violet, dotted] (.788675, 0) -- (.788675, 0.3);
	  \node[violet, below] at (.788675,0) { $u_3$};
  \draw[blue, dotted] (.908248, 0) -- (.908248, 0.3);
	  \node[blue, below] at (.908248,0) { $u_4$};
	
  \node[below left] at (0,0) {\small 0};
  \node[below] at (1,0) {\small 1};

  \node[blue, left] at (-0.025, 1/12) { $\sfrac{M}{3}$};
  \node[violet, left] at (-0.025, 1/6) { $\sfrac{2M}{3}$};

  \node[red] at (0.5,0.325) {$\eta_{\rm FKPP}$};
\end{tikzpicture}
	\begin{tikzpicture}[scale=5.5]

	\tikzmath
		{
			\ay = .0917517;
			\bee = .2113249;
			\dee = .816497;
			\em =  -\dee/(\bee-\ay);
			\el = - (6/(\bee-\ay)^3)*((1/36) - \dee*(\bee-\ay) - (\em/2)*(\bee-\ay)^2);
		} 
  \draw[->] (-0.0125,0) -- (1.1,0) node[below right] {$u$};
  \draw[->] (0,-0.0125) -- (0,0.3) node[left] {};

  \draw[red,thick, domain=0:\ay] plot (\x, {\x*(1-\x)});
  \draw[red,thick, domain=\ay:\bee]
  	plot (\x, {(1/12) + \dee*(\x-\ay) + (\em/2)*(\x-\ay)^2 + (\el/6)*(\x-\ay)^2*(2*\x - 3*\bee + \ay)});
	
	\tikzmath{
			\ayb = .908248;
			\beeb = .788675;
			\deeb = -.816497;
			\emb =  -\deeb/(\beeb-\ayb);
			\elb = - (6/(\beeb-\ayb)^3)*((1/36) - \deeb*(\beeb-\ayb) - (\emb/2)*(\beeb-\ayb)^2);
		}
  \draw[red,thick, domain= \bee:\beeb] plot(\x, 1/9);

  \draw[red,thick, domain=\beeb:\ayb]
  	plot (\x, {(1/12) + \deeb*(\x-\ayb) + (\emb/2)*(\x-\ayb)^2 + (\elb/6)*(\x-\ayb)^2*(2*\x - 3*\beeb + \ayb)});
	
  \draw[red,thick, domain=.90824:1] plot (\x, {\x*(1-\x)});

  \draw[blue, dotted] (-0.025, 1/12) -- (1.1, 1/12);
  \draw[violet, dotted] (-0.025, 1/6) -- (1.1, 1/6);
  \draw[red, dotted] (-0.025,1/9) -- (1.1, 1/9)
  	node [right] {$\sfrac{4M}{9}$};

  \draw[blue, dotted] (.0917517, 0) -- (.09175, 0.3);
	  \node[blue, below] at (.0917517,0) { $u_1$};
  \draw[violet, dotted] (.2113249, 0) -- (.2113249, 0.3);
	  \node[violet, below] at (.2113249,0) { $u_2$};
  \draw[violet, dotted] (.788675, 0) -- (.788675, 0.3);
	  \node[violet, below] at (.788675,0) { $u_3$};
  \draw[blue, dotted] (.908248, 0) -- (.908248, 0.3);
	  \node[blue, below] at (.908248,0) { $u_4$};
	  	
  \node[below left] at (0,0) {\small 0};
  \node[below] at (1,0) {\small 1};

  \node[blue, left] at (-0.025, 1/12) { $\sfrac{M}{3}$};
  \node[violet, left] at (-0.025, 1/6) { $\sfrac{2M}{3}$};
  
  \node[red] at (0.5,0.325) {$\psi\circ \eta_{\rm FKPP}$};
\end{tikzpicture}
\caption{An illustration of the definition of $u_1$, $u_2$, $u_3$, and $u_4$ as well as the shape of $\psi \circ \eta_{\rm FKPP}$.}
\label{f.eta}
\end{figure}

	Let $\psi$ 
	be a smooth, nondecreasing, concave function such that
	\be
		\psi(x) = \begin{cases}
				x
					\qquad &\text{ if } x \leq \frac{M}{3},\\
				\sfrac{4M}{9}
					\qquad &\text{ if } x \geq \frac{2M}{3}.
			\end{cases}
	\ee
	Then, note that 
	\be\label{e.c060602}
		\psi\circ \eta_{\rm FKPP}\in C^3\Big(\Big[\frac{u_1}{2}, \frac{1+u_4}{2}\Big]\Big)
	\ee
and
	\be
		\psi(\eta_{\rm FKPP}(u)) = \begin{cases}
				\eta_{\rm FKPP}(u)
					\qquad &\text{ if } u \in [0,u_1]\cup [u_4,1],\\
				\frac{4M}{9}
					\qquad &\text{ if } u \in [u_2,u_3].
			\end{cases}
	\ee
	See \Cref{f.eta} for an illustration of this.  
	
Let $\phi$ be a nonnegative, $C^3$-function such that  $\phi(u)>0$ for $u\in(u_2,u_3)$ and vanishes
outside of this interval, and 
	\be\label{e.c060601}
		\begin{split}
			&
			\quad
			\phi'\Big(\frac{u_2+u_3}{2}\Big) = 1,
			\quad
			\sup |\phi'| < 2, 
			\quad
			\phi''\Big(\frac{u_2+u_3}{2}\Big) > 0.
		\end{split}
	\ee
	We note that
	\be\label{e.c060603}
		\supp\big( ( \psi \circ \eta_{\rm FKPP})'\big)
			\cap \supp(\phi')
			= \emptyset.
	\ee
	
	We now define the nonlinearity.  Inspired by~\eqref{e.f_eta}, we let
	\be
		f(u) = \eta(u) (2 - \eta'(u))
	\ee
	where
	\be
		\eta(u) = \psi(\eta_{\rm FKPP}(u)) + \phi(u).
	\ee
	
	Let us first check that $f\in C^2([0,1])$.  This is clear away from $u=0$ and $u=1$. 
	On the other hand, for $u \in [0,u_1] \cup [u_4,1]$,
	\be
		\eta(u)
			= \eta_{\rm FKPP}(u)
	\ee
	and $\phi(u)=0$, hence, for all $u \in [0,u_1]\cup[u_4,1]$,
	\be\label{e.c060702}
		f(u) = \eta(u) (2 - \eta'(u))
			=  \eta_{\rm FKPP}(u) (2 - \eta_{\rm FKPP}'(u))
			= f_{\rm FKPP}(u) = u(1-u).
	\ee
	Thus, $f\in C^2([0,1])$. It also follows that $f(0) = f(1) = 0$ and $f'(0) = 1$, so that conditions 
We also note that, by construction, $\eta > 0$ and $\eta' < 2$, so that $f>0$.  Here we used the concavity of $\psi$, (\ref{e.c060601}),
and~\eqref{e.c060603} to guarantee that $\eta' < 2$.
	
	We note that, by construction, $\eta$ is the traveling wave profile function for $f$.  Indeed, defining $U$ to be the solution of
	\be
		- U' = \eta(U)
		\quad\text{ such that } U(0) = \frac12,
	\ee
	then $U$ is a speed two traveling wave solution of~\eqref{e.rde}.  As $\eta(u)>0$ for $u\in(0,1)$ and $\eta(0)=\eta(1)=0$, 
	we also know that $U(-\infty) = 1$ and $U(+\infty) = 0$.  Additionally, since $f'(0) = 1$,  the minimal speed is~$c_*=2$.  
	Hence, $U$ is the minimal speed traveling wave of~\eqref{e.rde}, and $\eta$ is the minimal speed traveling wave profile function associated to $f$.

	The proof is then complete by noting that $\eta' \equiv \phi'$ near $\sfrac{(u_2+u_3)}{2}$
	and, hence,
	\be
		\begin{split}
		Q\Big(\frac{u_2+u_3}{2}\Big)
			&= \phi'\Big(\frac{u_2+u_3}{2}\Big)\Big(2 - \phi'\Big(\frac{u_2+u_3}{2}\Big)\Big)
				+ \eta\Big(\frac{u_2+u_3}{2}\Big) \phi''\Big(\frac{u_2+u_3}{2}\Big)
			\\&
			= 1 + \phi\Big(\frac{u_2+u_3}{2}\Big) \phi''\Big(\frac{u_2+u_3}{2}\Big)
			> 1,
		\end{split}
	\ee
	where we used~\eqref{e.c060601} in the second equality and in the inequality.~$\Box$

\bibliographystyle{plain}

\end{document}